\documentclass[aos,MSNbibl,seceqn,dvips]{arximspdf}
\usepackage{mathbh,dcolumn,upgreek}
\usepackage{graphicx}

\doi{10.1214/12-AOS1010} 
\volume{40}
\issue{3}
\pubyear{2012}
\firstpage{1376}
\lastpage{1402}

\makeatletter
\newcolumntype{d}[1]{D{.}{.}{#1}}

\newcommand{\rrVert}{\Vert}
\newcommand{\rrvert}{\vert}
\newcommand{\llVert}{\Vert}
\newcommand{\llvert}{\vert}
\def\II{\mathit{I I}}
\def\III{\mathit{I I I}}
\newtheorem{lemma}{Lemma}[section]
\newtheorem{proposition}{Proposition}[section]
\newtheorem{theorem}{Theorem}[section]
\makeatother

\begin{document}
\begin{frontmatter}

\title{Inference of time-varying regression models\thanksref{T1}}
\thankstext{T1}{Supported in part by NSF Grant DMS-09-06073.}
\runtitle{Inference of TVCM}

\begin{aug}
\author[A]{\fnms{Ting} \snm{Zhang}\corref{}\ead[label=e1]{ting-zhang@uiowa.edu}}
\and
\author[B]{\fnms{Wei Biao} \snm{Wu}\ead[label=e2]{wbwu@galton.uchicago.edu}}
\runauthor{T. Zhang and W. B. Wu}
\affiliation{University of Chicago}
\address[A]{Department of Statistics \\
\quad and Actuarial Science\\
University of Iowa\\
Iowa City, Iowa 52242-1409\\
USA\\
\printead{e1}}
\address[B]{Department of Statistics\\
University of Chicago\\
5734 S. University Ave.\\
Chicago, Illinois 60637\\
USA\\
\printead{e2}}
\end{aug}

\received{\smonth{2} \syear{2012}}
\revised{\smonth{4} \syear{2012}}

%
\begin{abstract}
We consider parameter estimation, hypothesis testing and variable
selection for partially time-varying coefficient models. Our asymptotic
theory has the useful feature that it can allow dependent,
nonstationary error and covariate processes. With a two-stage method,
the parametric component can be estimated with a $n^{1/2}$-convergence
rate. A simulation-assisted hypothesis testing procedure is proposed
for testing significance and parameter constancy. We further propose an
information criterion that can consistently \mbox{select} the true set of
significant predictors. Our method is applied to autoregressive models
with time-varying coefficients. Simulation results and a real data
application are provided.
\end{abstract}

%
\begin{keyword}[class=AMS]
\kwd[Primary ]{62G05}
\kwd{62G10}
\kwd[; secondary ]{62G20}.
\end{keyword}

\begin{keyword}
\kwd{Information criterion}
\kwd{locally stationary processes}
\kwd{nonparametric hypothesis testings}
\kwd{time-varying coefficient models}
\kwd{variable selection}.
\end{keyword}

\end{frontmatter}

\section{Introduction}\label{secintro}
Varying coefficient models have been extensively studied in the
literature, and they are useful for characterizing nonconstancy
relationship between predictors and responses in regression models;
see, for example, \cite
{ZegerDiggle1994,HooverRiceWuYang1998,FanZhang1999,FanZhang2000,LinYing2001,ZhangLeeSong2002,HuangWuZhou2004,RamsaySilverman2005}.
In this paper we consider the time-varying coefficient model
%
\begin{equation}
\label{eqnlreg} y_i = \mathbf{x}_i^\top
\bolds{\beta}_i + e_i,\qquad i = 1,\ldots,n,
\end{equation}
where $y_i$ is the response, $\mathbf{x}_i$ is the predictor, $^\top
$ is the transpose operator, $\bolds{\beta}_i = \bolds{\beta}
(i/n)$ for some smooth function $\bolds{\beta}\dvtx  [0,1] \to\mathbb
R^p$ and $e_i$ is the error. We assume that $E(e_i |\mathbf
{x}_i) = 0$. Estimation of the coefficient function $\bolds{\beta}
(\cdot)$ in model~(\ref{eqnlreg}) has been considered by \cite
{Robinson1989,Robinson1991,OrbeFerreiraRodriguez2005,Cai2007,ZhouWu2010}
among others. An important special example of (\ref{eqnlreg}) is the
time-varying autoregressive model \cite
{DahlhausNeumannSachs1999,MoulinesPriouretRoueff2005} by letting
$\mathbf{x}_i = (y_{i-1},\ldots,y_{i-p})^\top$. There are
important differences between our model (\ref{eqnlreg}) and those
under longitudinal or functional setting. Here we assume that only one
realization $(\mathbf{x}_i, y_i)_{i=1}^n$ is available, while in the
longitudinal and functional setting, many subjects are measured at
multiple times, so that one has multiple realizations.

A natural problem for model (\ref{eqnlreg}) is to test whether
certain or all components in $\bolds{\beta}_i$ are time-invariant.
There is a huge literature on the problem of testing parameter
stability; see, for example, \cite
{Chow1960,BrownDurbinEvans1975,NabeyaTanaka1988,LeybourneMcCabe1989,Nyblom1989,PlobergerKramerKontrus1989,Andrews1993,LinTerasvirta1994,DavisHuangYao1995,FanZhangZhang2001,ZhangDette2004,FanHuang2005,HeTerasvirtaGonzalez2009}.
For model~(\ref{eqnlreg}), we are interested in testing
%
\begin{equation}
\label{eqnH0} H_0\dvtx  \mathbf{A} \bolds{\beta}(\cdot) \equiv
\mathbf{a}
\end{equation}
for some vector $\mathbf{a} \in\mathbb{R}^s$, where $\mathbf{A}
\in\mathbb{R}^{s \times p}$ is a real-valued matrix. With an
appropriately chosen $\mathbf{A}$, the null hypothesis (\ref
{eqnH0}) can be formulated to test whether a certain part of
coefficients is zero or time-invariant. In the latter case,
$\mathbf{a}$ needs to be replaced by an estimate $\hat{\mathbf{a}}$. Zhou and Wu \cite{ZhouWu2010} built simultaneous confidence
tubes for the regression coefficient function $\bolds{\beta}(\cdot
)$, which can be used as an $\mathcal L^\infty$-test for (\ref
{eqnH0}). The latter test often does not have a good power if the
alternative hypothesis consists of nonzero smooth functions. In Section
\ref{subsechypotest} we propose a more powerful $\mathcal L^2$-test
which is based on the weighted integrated squared errors. Our setting
is much more general than the one in \cite{ChenHong2012} where
$(\mathbf{x}_i, e_i)$ is assumed to be $\beta$-mixing and
stationary. In comparison, we allow nonstationary predictor and error
processes which can be nonstrong mixing; see Section \ref{secprelim}
for our nonstationary framework and basic assumptions.

If some of the coefficients are time-invariant, model (\ref{eqnlreg})
becomes the (semiparametric) partially time-varying coefficient model
%
\begin{equation}
\label{eqnregsemi} y_i = \mathbf{x}_{D_1,i}^\top
\bolds{\beta}_{D_1} + \mathbf{x}_{D_2,i}^\top\bolds{
\beta}_{D_2}(i/n) + e_i, \qquad i = 1,\ldots,n,
\end{equation}
where $D_1, D_2 \subseteq D^* = \{1,\ldots,p\}$ are groups of
parametric and nonparametric components, respectively. Based on an
estimate of (\ref{eqnlreg}), we simply take an integration/average
over the parametric part to obtain an estimate of~$\bolds{\beta}_{D_1}$ that achieves the $n^{1/2}$-convergence rate. An asymptotic
theory is given in Section \ref{subsecparameterestimation}. This
method was
previously used in \cite{ZhangLeeSong2002} for estimating
state-domain semivarying coefficient models. The latter paper assumed
that $y_i = \mathbf{x}_i^\top\bolds{\beta}(\mathbf{u}_i) +
e_i$ where $(y_i, \mathbf{x}_i, \mathbf{u}_i)$, $i = 1,\ldots
,n$, are independent and identically distributed; see \cite
{XiaZhangTong2004} for the case with stationary mixing processes.
Our time-domain model (\ref{eqnregsemi}) is very general, and it
includes both (\ref{eqnlreg}) and usual linear regression models. Gao
and Hawthorne \cite{GaoHawthorne2006} considered a special example
of (\ref{eqnregsemi}) with $D_1=\{2, \ldots, p\}$, $D_2 = 1$ and
$\mathbf{x}_{D_2,i} \equiv1$, so that only the intercept term in
(\ref{eqnregsemi}) is time-varying.

Section \ref{subsecvarsel} deals with the problem of selecting
significant predictors. Fan et al. \cite{FanYaoCai2003} proposed an
extended AIC for choosing locally significant variables. Abramovich et
al. \cite{AbramovichSpencerTurley2007} considered the problem of
order selection for time-varying autoregressive models by requiring
that multiple realizations are available. Using the dependence
framework in Section \ref{secprelim}, we are able to solve this
problem under parameter instability and temporal dependence. In
particular, we propose an information criterion, consisting of measures
of goodness-of-fit and model complexity, that can consistently select
the true set of relevant predictors based only on one realization.
Section \ref{secimplementation} provides simulation studies and an
application. Proofs are given in the \hyperref[secproof]{Appendix}.

\section{Model assumptions}\label{secprelim}
Since the coefficient function $\bolds{\beta}(\cdot)$ in (\ref
{eqnlreg}) is smooth, we can naturally estimate it (along with its
derivative) by
%
\begin{equation}
\label{eqnllmineqn} \bigl(\tilde{\bolds{\beta}}(t), \tilde{\bolds{
\beta}}'(t)\bigr) = \mathop{\operatorname{argmin}}_{\bolds{\eta}_0, \bolds{\eta}_1
\in\mathbb{R}^p} \sum
_{i=1}^n \bigl\{y_i -
\mathbf{x}_i^\top \bolds{\eta}_0 -
\mathbf{x}_i^\top\bolds{\eta}_1 (i/n - t)\bigr
\}^2 K \biggl(\frac{i/n - t}{b_n} \biggr),\hspace*{-35pt}
\end{equation}
where $K(\cdot)$ is the kernel function, and $b_n$ is a bandwidth
sequence satisfying $b_n \to0$ and $nb_n \to\infty$. Throughout the
paper we assume that the kernel function $K(\cdot)$ is a symmetric and
bounded function in $\mathcal C^1[-1,1]$ with $\int_{-1}^1 K(v) \,dv =
1$. For example, it can be the Epanechnikov kernel $K(v) = 3 \max(0,
1-v^2)/4$ or the Bartlett kernel $K(v) = \max(0, 1-|v|)$. Observe that
(\ref{eqnllmineqn}) has the closed form solution
%
\begin{equation}\qquad
\label{eqnllest} \pmatrix{ \tilde{\bolds{
\beta}}(t)
\vspace*{2pt}\cr
b_n\tilde{\bolds{\beta}}'(t) }=  \pmatrix{
 \mathbf{U}_{n,0}(t)
& \mathbf{U}_{n,1}(t)
\vspace*{2pt}\cr
\mathbf{U}_{n,1}(t) & \mathbf{U}_{n,2}(t) }^{-1} \pmatrix{
\mathbf{V}_{n,0}(t)
\vspace*{2pt}\cr
\mathbf{V}_{n,1}(t) } =
\mathbf{U}_n(t)^{-1} \mathbf{V}_n(t),
\end{equation}
where for $l \in\{0,1,2\}$,
\[
\mathbf{U}_{n,l}(t) = (nb_n)^{-1} \sum
_{i=1}^n \mathbf{x}_i
\mathbf{x}_i^\top\bigl\{(i/n - t)/b_n\bigr
\}^l K\bigl\{(i/n-t)/b_n\bigr\},
\]
with the convention that $0^0 = 1$, and
\[
\mathbf{V}_{n,l}(t) = (nb_n)^{-1} \sum
_{i=1}^n \mathbf{x}_i y_i
\bigl\{(i/n - t)/b_n\bigr\}^l K\bigl\{(i/n-t)/b_n
\bigr\}.
\]
To establish an asymptotic theory for $\tilde{\bolds{\beta}
}(\cdot)$, we need to impose appropriate regularity conditions on the
covariates $(\mathbf{x}_i)$ and errors $(e_i)$. For testing the
hypothesis~(\ref{eqnH0}), \cite{ChenHong2012} assumed that
$(\mathbf{x}_i, e_i)$ is $\beta$-mixing and stationary. To allow
nonstationary predictor and error processes that can be nonstrong
mixing, we assume that
%
\begin{equation}
\label{eqnformulation} \mathbf{x}_i = \mathbf{G}(i/n; \bolds{
\mathcal{F}}_i)\quad \mathrm{and}\quad e_i = H(i/n; \bolds{
\mathcal{F}}_i),
\end{equation}
where $\bolds{\mathcal{F}}_i = (\ldots,\bolds{\varepsilon}_{i-1},\bolds{\varepsilon}_i)$ is a shift process of independent and
identically distributed (i.i.d.) random variables $\bolds{\varepsilon}_k$, $k \in\mathbb{Z}$ and $\mathbf{G}$ and $H$ are measurable
functions such that $\mathbf{G}(t;\bolds{\mathcal{F}}_i)$ and
$H(t;\bolds{\mathcal{F}}_i)$ are well defined for each $t \in
[0,1]$. This setup is also used in \cite{ZhouWu2010}.

For a random vector $\mathbf Z$, we write $\mathbf Z \in
\mathcal L^q$, $q > 0$, if $\|\mathbf Z\|_q = \{E(|\mathbf
Z|^q)\}^{1/q} < \infty$, where $|\cdot|$ is the Euclidean vector
norm, and we denote $\|\cdot\| = \|\cdot\|_2$. A~process $\mathbf{J}(t; \bolds{\mathcal{F}}_k)$ is said to be stochastically
Lipschitz continuous ($\mathbf{J} \in \operatorname{Lip}$ in short) if there exists
$C > 0$, such that $\|\mathbf{J} (t_1; \bolds{\mathcal{F}}_k)
- \mathbf{J} (t_2; \bolds{\mathcal{F}}_k)\| \leq C |t_1 -
t_2|$ holds uniformly for all $t_1, t_2 \in[0,1]$. Then, under
condition (A2) below,~(\ref{eqnformulation}) defines locally
stationary processes. Let $\{\bolds{\varepsilon}_j'\}_{j \in\mathbb
Z}$ be an i.i.d. copy of $\{\bolds{\varepsilon}_j\}_{j \in\mathbb
Z}$ and $\bolds{\mathcal F}_{i,\{0\}} = (\ldots, \bolds{
\varepsilon}_{-1}, \bolds{\varepsilon}_0', \bolds{\varepsilon}_1,
\ldots, \bolds{\varepsilon}_i)$ be the coupled shift process. We
define the functional dependence measure
\[
\delta_{k,q}(\mathbf{J}) = \sup_{t \in[0,1]} \bigl\|\mathbf{J}(t; \bolds{
\mathcal F}_k) - \mathbf{J} (t; \bolds{\mathcal F}_{k,\{0\}})
\bigr\|_q \quad\mathrm{and}\quad \Theta_{m,q}(\mathbf{J}) = \sum
_{j=m}^{\infty} \delta_{j,q}(\mathbf{J}).
\]
Let $\bolds{\Lambda}(\mathbf{J}, t) = \sum_{k \in\mathbb Z}
\operatorname{cov}\{\mathbf{J} (t; \bolds{\mathcal{F}}_0),
\mathbf{J} (t; \bolds{\mathcal{F}}_k)\}$ be the long-run
covariance matrix, and $\mathbf{M}(\mathbf{J}, t) = E\{
\mathbf{J}(t; \bolds{\mathcal{F}}_0) \mathbf{J}(t;
\bolds{\mathcal{F}}_0)^\top\}$. Under the short-range
dependence condition $\Theta_{0,2}(\mathbf{J}) < \infty$, both of
them are uniformly bounded over $t \in[0,1]$. Let $\mathbf{L}(t;
\bolds{\mathcal{F}}_k) = \mathbf{G}(t; \bolds{\mathcal
{F}}_k) H(t; \bolds{\mathcal{F}}_k)$, and we shall make the
following assumptions:
\begin{longlist}[(A1)]
\item[(A1)] Smoothness: $\bolds{\beta}\in\mathcal{C}^3[0,1]$;
\item[(A2)] Local stationarity: $\mathbf{G}, \mathbf{L} \in \operatorname{Lip}$;
\item[(A3)] Short-range dependence: $\Theta_{0,4}(\mathbf{G}) +
\Theta_{0,\iota}(\mathbf{L}) < \infty$ for some $\iota> 2$;
\item[(A4)] The smallest eigenvalue of $\mathbf{M}(\mathbf{G},
\cdot)$ is bounded away from zero on $[0,1]$.
\end{longlist}
A sufficient condition for (A3) is that $\Theta_{0,2\iota
}(\mathbf{G}) + \Theta_{0,2\iota}(H) < \infty$ for some $\iota> 2$.

\section{Main results}\label{secmain}
\subsection{Parameter estimation}\label{subsecparameterestimation}
Let $\mathbf{A}$ be a pre-specified matrix and $\mathbf{a} =  \int_0^1 \mathbf{A} \bolds{\beta}(t) \,dt$. Then
\[
\hat{\mathbf{a}} = \int_0^1 \mathbf{A} \tilde{
\bolds \beta}(t) \,dt
\]
is an estimate of $\mathbf{a}$. For the partially time-varying
coefficient model (\ref{eqnregsemi}), let $\mathbf{A}_{D_1} \in
\mathbb{R}^{p_1 \times p}$ be a matrix with rows $\{\mathbf
z_i^\top\}_{i \in D_1}$, where $\mathbf z_i \in\mathbb{R}^p$ is
the vector with unit $i$th component and zeros elsewhere, then
$\mathbf{A}_{D_1} \bolds{\beta}(t) = \bolds{\beta}_{D_1}(t)$, $t \in[0,1]$. Although $\bolds{\beta}_{D_1}$ can be
consistently estimated by $\tilde{\bolds{\beta}}_{D_1}(t)$ for
any $t \in(0,1)$, the smoothed estimate
\[
\hat{\bolds{\beta}}_{D_1} = \int_0^1
\mathbf{A}_{D_1} \tilde {\bolds{\beta}}(t) \,d t
\]
can have a better rate of convergence.

%
\begin{theorem}\label{thmCLTahat}
Assume \textup{(A1)--(A4)} and $\Theta_{n,\iota}(\mathbf{L}) =
\mathcal{O}(n^{-\nu})$ for some $\nu> 1/2-1/\iota$. Let
$\bolds\Xi(t) = \mathbf{M}(\mathbf{G},t)^{-1}
\bolds{\Lambda}(\mathbf{L},t) \mathbf{M}(\mathbf{G},t)^{-1}$ and\vadjust{\goodbreak} $\kappa_2 = \int_{-1}^1 v^2 K(v) \,d v$. If $b_n
\asymp n^{-c}$ for some $1/6 < c < \min\{1/3,1/2-1/(2\iota)\}$,
then
\[
n^{1/2}(\hat{\mathbf{a}} - \mathbf{a} - \bolds\xi_n)
\Rightarrow N \biggl\{0, \int_0^1 \mathbf{A}
\bolds\Xi(t) \mathbf{A}^\top \,dt \biggr\} \quad\mbox{where } \bolds
\xi_n = \frac{b_n^2 \kappa_2}{2} \int_0^1
\mathbf{A} \bolds \beta''(t) \,dt.
\]
\end{theorem}

In Theorem \ref{thmCLTahat}, the term $\bolds\xi_n$ can be
interpreted as the bias due to nonparametric estimation, and it
vanishes under the null hypothesis (\ref{eqnH0}). Hence the
parametric component ${\bolds{\beta}}_{D_1}$ in the
semi-parametric model (\ref{eqnregsemi}) can have a
$n^{1/2}$-consistent estimate $\hat{\bolds{\beta}}_{D_1}$.

\subsection{Hypothesis testing}\label{subsechypotest}
For testing the null hypothesis (\ref{eqnH0}), let $\mathbf{W}(\cdot)$ be
a continuous mapping from $[0,1]$ to symmetric
positive-definite matrices in $\mathbb{R}^{s \times s}$. Consider the
weighted integrated squared error
%
\begin{equation}
\label{eqnTn} T_n(\mathbf{A},\mathbf{a},\mathbf{W}) = \int
_0^1 \bigl\{ \mathbf{A} \tilde{\bolds{\beta}}(t)
- \mathbf{a}\bigr\}^\top \mathbf{W}(t) \bigl\{\mathbf{A} \tilde{\bolds{
\beta}}(t) - \mathbf{a}\bigr\} \,dt.
\end{equation}
If $\mathbf{a}$ is unknown, an estimate can be used. For example we
can use $\hat{\mathbf{a}} = \int_0^1 \mathbf{A} \tilde
{\bolds{\beta}}(t) \,dt$, which has the parametric convergence rate;
see Theorem \ref{thmCLTahat}. Chen and Hong \cite{ChenHong2012}
considered the special case that $(\mathbf{x}_i, e_i)$ is a
stationary $\beta$-mixing process. Their generalized Hausman test
\cite{Hausman1978} relates to (\ref{eqnTn}) with $\mathbf{A}$
being the identity matrix and $\mathbf{W}(t) = \mathbf{M}(\mathbf{G},t)$. Such a choice of weight matrices should be used if
we are interested in prediction. Alternatively, we can use $\mathbf{W}(t) \equiv\mathbf I_{s \times s}$, the identity matrix to form
the integrated squared errors. Let $K_2 = \int_{-1}^1 K(v)^2 \,dv$, by
Theorem~1 in \cite{ZhouWu2010}, $\mathbf{A} \tilde{\bolds
\beta}(t)$ has the asymptotic covariance $(nb_n)^{-1} K_2 \mathbf{A} \bolds\Xi(t) \mathbf{A}^\top$. Hence, we can choose
$\mathbf{W}(t) = \{\mathbf{A} \bolds\Xi(t) \mathbf{A}^\top\}^{-1}$ to serve as a normalizer. In this case, (\ref{eqnTn})
is (proportionally) an integral of the squared local $t$-statistics.

For a matrix $\mathbf{A}$, define $\underline\rho(\mathbf{A}) = \inf\{|\mathbf{A} \mathbf{v}|\dvtx  |\mathbf{v}| = 1
\}$ and $\overline\rho(\mathbf{A}) = \sup\{|\mathbf{A}
\mathbf{v}|\dvtx \allowbreak |\mathbf{v}| = 1 \}$. Let
\[
K^*(x) = \int_{-1}^{1-2|x|} K(v) K\bigl(v+2|x|\bigr) \,dv,
\]
and $K^*_2 = \int_{-1}^1 K^*(v)^2 \,dv$. Since $K \in\mathcal{K}$,
we have $K^* \in\mathcal{C}^1[-1, 1]$ and is symmetric. Let
%
\begin{eqnarray}
\label{thmO19303} \bolds\Xi_{\mathbf{A},\mathbf{W}}(t) &=& \mathbf{W}(t)^{1/2}
\mathbf{A} \bolds\Xi(t) \mathbf{A}^\top \mathbf{W}(t)^{1/2},
\nonumber
\\[-8pt]
\\[-8pt]
\nonumber
\Xi_{\mathbf{A},\mathbf{W},l}& =& \operatorname{tr} \biggl\{\int_0^1
\bolds\Xi_{\mathbf{A},\mathbf{W}}(t)^l \,dt \biggr\}.
\end{eqnarray}
Theorem \ref{thmCLTTn} provides asymptotic normality for
$T_n(\mathbf{A},\mathbf{a},\mathbf{W})$.

%
\begin{theorem}\label{thmCLTTn}
Assume \textup{(A1)--(A4)}, $\Theta_{0,4}(\mathbf{L}) < \infty$ and
$\Theta_{n,\iota}(\mathbf{L}) = \mathcal{O}(n^{-\nu})$ for some
$\nu> 1$. If $b_n \asymp n^{-c}$ for some\vadjust{\goodbreak} $2/11 < c < \min\{1/3,
3/5-4/(5\iota),\break 2-4/\iota\}$, then
%
\begin{equation}\qquad
\label{eqnCLTTn} nb_n^{1/2} \bigl\{T_n(
\mathbf{A},\mathbf{a},\mathbf{W}) - (nb_n)^{-1} K^*(0)
\Xi_{A,\mathbf{W},1}\bigr\} \Rightarrow N \bigl(0, 4 K^*_2
\Xi_{\mathbf{A},\mathbf{W},2}\bigr).
\end{equation}
If in addition $\hat{\mathbf{a}} = \mathbf{a} +
O_p(n^{-1/2})$, then (\ref{eqnCLTTn}) holds for $T_n(\mathbf{A},\hat{\mathbf{a}},\mathbf{W})$.
\end{theorem}

Let $\Phi(\cdot)$ be the cumulative standard normal distribution function
and $q_{1-\alpha}$ be the corresponding
$(1-\alpha)$th quantile. We reject the null hypothesis~(\ref{eqnH0}) at level $\alpha$ if
%
\begin{equation}
\label{eqntestusual}  T_n(\mathbf{A},\hat{\mathbf{a}},\mathbf{W}) >
b_n^{-1/2} K^*(0) \Xi_{\mathbf{A},\mathbf{W},1} + n^{-1}
b_n^{-1/2} \bigl(4 K^*_2 \Xi_{\mathbf{A},\mathbf{W},2}
\bigr)^{1/2} q_{1-\alpha}.\hspace*{-35pt}
\end{equation}
Let $\mathbf{f}\dvtx  [0,1] \to\mathbb{R}^{s}$ be of class
$\mathcal C^3$, and $\{d_n\}$ be a sequence of nonnegative real
numbers. Proposition \ref{thmPower} provides the asymptotic power
of the test (\ref{eqntestusual}) under the local alternative
%
\begin{equation}
\label{eqnlocalalternative} \mathbf{A} \bolds{\beta}(t) = \mathbf{a} +
d_n \mathbf{f}(t).
\end{equation}

%
\begin{proposition}\label{thmPower}
Assume conditions of Theorem \ref{thmCLTTn}. If $nb_n^{1/2} \,d_n^2
\to\break s > 0$, then the power of the test (\ref{eqntestusual})
satisfies
%
\begin{equation}
\label{eqnPower} \mathrm{Power} \to\Phi \biggl\{q_\alpha+
\frac{s \int_0^1 \mathbf{f}(t)^\top\mathbf{W}(t) \mathbf{f}(t) \,dt }{(4 K_2^*
\Xi_{\mathbf{A},\mathbf{W},2})^{1/2}} \biggr\}.
\end{equation}
\end{proposition}

\subsection{Variable selection}
\label{subsecvarsel} In this section we shall propose an
information criterion for time-varying coefficient models that can
consistently identify the true set of relevant predictors. Recall
that $D^* = \{1,\ldots,p\}$ is the whole set of potential
predictors, and $\tilde{\bolds{\beta}}(\cdot)$ is the
coefficient estimate. Let $D_0$ be the true set of relevant
predictors. For a candidate subset $D \subseteq D^*$, we can
compute the variable selection information criterion
%
\begin{eqnarray}
\label{eqnVIC} \operatorname{VIC}(D) = \log\bigl\{\operatorname{RSS}(D)\bigr\} +
\chi_n |D|
\nonumber
\\[-8pt]
\\[-8pt]
 \eqntext{\mbox{where } \displaystyle\operatorname{RSS}(D) = \sum
_{i=1}^n \bigl\{y_i -
\mathbf{x}_{D,i}^\top\tilde{\bolds{\beta}}_D(i/n)
\bigr\}^2.}
\end{eqnarray}
Here $\chi_n$ is a tuning parameter. We select a subset $D$ that
minimizes $\operatorname{VIC}(D)$, thus balancing goodess-of-fit and
model complexity. Smaller $\chi_n$ leads to more predictors, and
vice versa. Theorem \ref{thmmodelselect} provides theoretical
properties of our procedure.

%
\begin{theorem}\label{thmmodelselect}
Assume \textup{(A1)--(A4)}, $\Theta_{0,4}(\mathbf{L}) + \Theta_{0,2}(H)
< \infty$, $\Theta_{n,\iota}(\mathbf{L}) =
\mathcal{O}(n^{-\nu})$ for some $\nu> 1/2-1/\iota$. Let
$\varphi_n = (nb_n)^{-1} \{n^{1/\iota} + (nb_n \log n)^{1/2}\} +
b_n^2$ and $\rho_n = n^{-1/2} b_n^{-1} + b_n$. If $b_n \asymp
n^{-c}$ for some $0 < c < 1-1/\iota$,
%
\begin{equation}
\label{eqnchin} \chi_n \to0 \quad\mbox{and}\quad \bigl\{\varphi_n(
\varphi_n + \rho_n)\bigr\}^{-1}
\chi_n \to\infty,
\end{equation}
then, for any $D \neq D_0$, $\operatorname{pr} \{\operatorname{VIC}(D) >
\operatorname{VIC}(D_0) \} \to1$.
\end{theorem}

\section{Implementation}\label{secimplementation}
\subsection{Covariance matrix estimation}\label{subseccovest}
$\!\!\!$Theorems \ref{thmCLTahat} and \ref{thmCLTTn} both involve~un\-known
quantities depending on the covariance matrices $\mathbf{M}(\mathbf{G},t)$ and $\bolds{\Lambda}(\mathbf{L},t)$, $t \in
[0,1]$. The problem of estimating covariance matrices has been
extensively studied; see \cite
{NeweyWest1987,AndrewsMonahan1992,LumleyHeagerty1999} among
others. Let $\varpi_n$, $\tau_n$ and $\varrho_n$ be bandwidth
sequences satisfying $\varpi_n \to0$, $\tau_n \to0$, $\varrho_n
\to0$ and $n \tau_n \varrho_n \to\infty$. Let $\mathcal{I}_{n,1}
= [0,\tau_n \varrho_n]$, $\mathcal{I}_{n,2} = (\tau_n \varrho_n, 1
- \tau_n \varrho_n)$, $\mathcal{I}_{n,3} = [1 - \tau_n \varrho_n,1]$ and
%
\[
\bolds{\lambda}_i(\mathbf{L},\tau_n
\varrho_n)
= \cases{ \displaystyle\mathbf{L}_i \mathbf{L}_i^\top+
2 \mathbf{L}_i \sum_{j=1}^n
\mathbf{L}_j^\top\mathbh{1}_{ \{ {0 < j/n - i/n \leq\tau_n
\varrho _n}  \} }, &\quad  $\mbox{if $i/n \in\mathcal{I}_{n,1}$;}$\vspace*{2pt}
\cr
\displaystyle\mathbf{L}_i
\sum_{j=1}^n \mathbf{L}_j^\top
\mathbh{1}_{
\{ {|i/n - j/n| \leq\tau_n \varrho_n}  \} }, &\quad  $\mbox{if $i/n \in \mathcal{I}_{n,2}$;}$
\vspace*{2pt}
\cr
\displaystyle\mathbf{L}_i \mathbf{L}_i^\top+
2 \mathbf{L}_i \sum_{j=1}^n
\mathbf{L}_j^\top\mathbh{1}_{ \{ {0 < i/n - j/n \leq\tau_n
\varrho _n}  \} }, &\quad  $\mbox{if
$i/n \in\mathcal{I}_{n,3}$.}$ }
\]
For $t \in[0,1]$, we estimate $\mathbf{M}(\mathbf{G},t)$ and
$\bolds{\Lambda}(\mathbf{L},t)$,
respectively, by
\[
\hat{\mathbf{M}}(\mathbf{G}, t) = \sum_{i = 1}^n
\mathbf{x}_i \mathbf{x}_i^\top
\omega_{i,\varpi_n}(t)
\]
and
\[
\hat{\bolds{\Lambda}}(\mathbf{L}, t) = \sum_{i = 1}^n
\frac{\bolds{\lambda}_i(\mathbf{L},\tau_n
\varrho_n) + \bolds{\lambda}_i(\mathbf{L},\tau_n
\varrho_n)^\top}{2} \omega_{i,\tau_n}(t),
\]
where $\omega_{i,b}(t) = K\{(i/n-t)/b\}
\{P_{b,2}(t)-(t-i/n)P_{b,1}(t)\}/\{P_{b,2}(t)P_{b,0}(t)-P_{b,1}(t)^2\}$
are local linear weights with bandwidth $b$ and $P_{b,l}(t) =
\sum_{j=1}^n (t-j/n)^l K\{(j/n-t)/b\}$. Proposition
\ref{thmcovest} provides consistency of our covariance
matrix estimates.

%
\begin{proposition}
\label{thmcovest} Assume \textup{(A2)}, $\Theta_{0,4}(\mathbf{G}) +
\Theta_{0,4}(\mathbf{L}) < \infty$ and
$\Theta_{n,2}(\mathbf{L}) = O(n^{-\nu})$ for some $\nu> 0$. If
both $\mathbf{M}(\mathbf{G},\cdot)$ and $\bolds
\Lambda(\mathbf{L},\cdot)$ are in class $\mathcal C^2$, then
%
\begin{equation}
\label{eqnbndM} \sup_{t \in[0,1]} \bigl\| \hat{\mathbf{M}}(\mathbf{G}, t) -
\mathbf{M}(\mathbf{G}, t) \bigr\| = O\bigl\{(n\varpi_n)^{-1/2} +
\varpi_n^2\bigr\},
\end{equation}
and
%
\begin{equation}
\label{eqnbndLambda} \sup_{t \in[0,1]} \bigl\| \hat{\bolds{\Lambda}}(\mathbf{L}, t)
- \bolds{\Lambda}(\mathbf{L}, t) \bigr\| = O\bigl\{\varrho_n^{1/2}
+ (n \tau_n \varrho_n)^{-\nu} + (
\tau_n \varrho_n)^{\nu/(1+\nu)} + \tau_n^2
\bigr\}.\hspace*{-35pt}
\end{equation}
\end{proposition}

The bound in (\ref{eqnbndM}) is optimized and becomes
$O(n^{-2/5})$ if $\varpi_n \asymp n^{-1/5}$. The optimal bound in
(\ref{eqnbndLambda}) is complicated and it depends on $\nu$,
the decay rate of dependence. In particular, if $\nu\geq
2/3$, then the optimal bound in (\ref{eqnbndLambda}) is
$O\{n^{-2\nu/(5\nu+2)}\}$ if $\tau_n \asymp
n^{-\nu/(5\nu+2)}$ and $\varrho_n \asymp n^{-4\nu/(5\nu+2)}$;
otherwise it is
$O\{n^{-\nu/(\nu+2)}\}$ if $\tau_n \asymp n^{-(1-\nu)/(\nu+2)}$
and $\varrho_n \asymp n^{-2\nu/(\nu+2)}$, or $\tau_n \asymp
n^{-\nu/(2\nu+4)}$ and $\varrho_n \asymp n^{-1/2}$. In computing
$\hat{\bolds{\Lambda}}({\mathbf{L}},t)$, since $(e_i)$
is usually unknown, we shall replace it by $(\tilde e_i)$,
the estimated local linear residuals.

\subsection{A simulation-assisted testing procedure}\label
{subsechypotestsimubased}
By the sandwich formula, let $\hat{\bolds\Xi}(t) =
\hat{\mathbf{M}}(\mathbf{G},t)^{-1} \hat{\bolds
\Lambda}(\tilde{\mathbf{L}},t) \hat{\mathbf{M}}(\mathbf{G},t)^{-1}$ and, as in (\ref{thmO19303}), correspondingly define
$\hat{\bolds\Xi}_{\mathbf{A},\mathbf{W}}(t)$ and $\hat
\Xi_{\mathbf{A},\mathbf{W},l}$. By (\ref{eqntestusual}), we
reject the null hypothesis (\ref{eqnH0}) at level $\alpha$ if
%
\begin{equation}
\label{eqO19307} \Delta_n(\mathbf{A},\hat{\mathbf{a}},\mathbf{W}) =
\frac{nb_n^{1/2} \{T_n(\mathbf{A},\hat{\mathbf{a}},\mathbf{W})
- (nb_n)^{-1} K^*(0) \hat\Xi_{\mathbf{A},\mathbf{W},1}\}
}{(4 K^*_2 \hat\Xi_{\mathbf{A},\mathbf{W},2})^{1/2}} > q_{1-\alpha}.\hspace*{-35pt}
\end{equation}
If $\mathbf{a}$ is known, then in (\ref{eqO19307}) we can use
$\mathbf{a}$ instead of $\hat{\mathbf{a}}$. The criterion
(\ref{eqO19307}) usually does not have a good performance because
of the slow convergence in (\ref{eqnCLTTn}). Note that the
statistic $\Delta_n(\mathbf{A}, \hat{\mathbf{a}}, \mathbf{W})$ is asymptotically pivotal, so we propose a
simulation-assisted testing procedure that can substantially improve
the finite-sample performance. In particular, we generate
i.i.d. standard normal random variables $y_i^\circ$, $i =
1,\ldots,n$, and i.i.d. standard multivariate normal random vectors
$\mathbf{x}_i^\circ$, $i = 1,\ldots,n$, that are also
independent of~$(y_k^\circ)$. We compute the corresponding $\Delta_n^\circ(\mathbf{A},\hat{\mathbf{a}},\mathbf{W})$, and
repeat this for many times to obtain its empirical quantile
$\hat q_{1-\alpha}$. We reject the null hypothesis~(\ref{eqnH0}) at
level $\alpha$ if $\Delta_n(\mathbf{A},\hat{\mathbf{a}},
\mathbf{W}) > \hat q_{1-\alpha}$. Our procedure has a similar
flavor as the Wilks type of phenomenon discussed in \cite
{FanZhangZhang2001}. A~major difference is that we allow dependent
and nonstationary errors.

\subsection{Bandwidth selection}\label{subsecbandwidth}
Bandwidth selection for nonparametric hypothesis testing is a
nontrivial problem, and it has been studied in \cite
{KulasekeraWang1997,FanLinton2003,DetteSpreckelsen2004,GaoGijbels2008}
among many others. As commented by Wang \cite{Wang2008}, there exists
no uniform guidance for an optimal choice. On
the positive side, our simulation results in Section
\ref{subsecsimulation} indicate that the empirical acceptance
probabilities are not quite sensitive to the choice of bandwidth.
Hence one can simply choose $b_n = n^{-1/5}$ that has the
asymptotic mean integrated squared error (AMISE) optimal rate. As
an alternative, we consider the generalized cross-validation (GCV)
selector by \cite{CravenWahba1979}, and estimate the covariance
matrix $\bolds\Gamma_n = \{E(e_ie_j)\}_{1 \leq i,j \leq n}$
to correct for dependence \cite{Wang1998}. Specifically, let
$\mathbf Y = (y_1,\ldots,y_n)^\top$; then for any bandwidth $b
\in(0,1)$, one can write the local linear smoothed fitted values
as $\hat{\mathbf Y}(b) = \mathbf{H}(b) \mathbf Y$,
where $\mathbf{H}(b)$ is the corresponding hat matrix. We
choose the bandwidth $\tilde b_n$ that minimizes
%
\begin{equation}
\label{eqnGCV} \operatorname{GCV}(b) = \frac{n^{-1} \{\hat{\mathbf Y}(b) - \mathbf
Y\}^\top\hat{\bolds\Gamma}{}^{-1}_n \{\hat{\mathbf Y}(b)
- \mathbf Y\} }{[1 - \operatorname{tr}\{\mathbf{H}(b)\}/n]^2}.
\end{equation}
An estimate of the covariance matrix $\bolds\Gamma_n$ can be
obtained by using the banding technique as in \cite
{BickelLevina2008,WuPourahmadi2009}. The GCV selector (\ref
{eqnGCV}) works reasonably well in our simulation studies.

We shall now provide data-driven choices of $\varpi_n$, $\tau_n$
and $\varrho_n$ in the estimation of covariance matrices. From the
construction in Section \ref{subseccovest}, we truncate the
long-run covariance matrix estimate at lag $m_n =
n\tau_n\varrho_n$ and, by the proof of Theorem~\ref{thmCLTahat},
\[
n^{1/2} \Biggl(\frac{1}{n} \sum_{i = 1}^{n-k}
\mathbf{L}_i^\top \mathbf{L}_{i+k} - \operatorname{tr}
\biggl[\int_0^1 \operatorname{cov}\bigl\{
\mathbf{L}(t; \bolds{\mathcal{F}}_0), \mathbf{L}(t; \bolds{
\mathcal{F}}_k)\bigr\}\,dt \biggr] \Biggr) \Rightarrow N\bigl(0,
\sigma_k^2\bigr),
\]
where $\sigma_k^2$ is the integrated long-run variance of the
process $\{\mathbf{L}_i^\top\mathbf{L}_{i+k}\}_{i=1}^{n-k}$. We propose to choose $\hat m_n = \max\{k
\geq0\dvtx  |n^{-1/2} \sum_{i = 1}^{n-k} \mathbf{L}_i^\top
\mathbf{L}_{i+k}| > 1.96 \sigma_k\}$. Note that the final
estimate is a local linear smoother of $\{\bolds
\lambda_i(\mathbf{L},\hat m_n/n) + \bolds
\lambda_i(\mathbf{L},\hat m_n/n)^\top\}/2$, $i = 1,\ldots,n$, and
we can apply the GCV method to select $\hat\tau_n$. The latter
can also be applied to $\mathbf{x}_i \mathbf{x}_i^\top$, $i
= 1,\ldots,n$, to select $\tilde\varpi_n$, and we take $\hat
\varpi_n = \max(\tilde\varpi_n,n^{-1/5})$ to avoid numerical
singularities. These data-driven choices of bandwidths are able to
capture dependence and nonstationarity and have a good
performance in our simulation studies.

For the information criterion (\ref{eqnVIC}), if $\iota\geq5/2$
and $b_n \asymp n^{-1/5}$, then (\ref{eqnchin}) becomes
%
\begin{equation}
\label{eqnchinsimplified} \chi_n \to0\quad \mathrm{and}\quad \bigl
\{n^{-3/5} (\log n)^{1/2}\bigr\}^{-1}
\chi_n \to\infty.
\end{equation}
Note that condition (\ref{eqnchinsimplified}) is more restrictive
than the traditional Bayesian information criterion (BIC) because
of parameter instability. Under the latter setting, a heavier
penalty on model complexity is usually needed to suppress the
over-fitting problem; see \cite{XiaZhangTong2004} for a similar finding
on cross-validation methods. As a rule of thumb, we suggest using
$\hat\chi_n = n^{-2/5}$. This simple choice performs reasonably
well as can be seen from our simulation study. The choice of
bandwidth becomes further complicated due to model uncertainty. We
suggest a two-stage selection procedure: let
$\tilde b_n$ be the selected bandwidth by GCV with all available
predictors, and we use the information criterion (\ref{eqnVIC})
to select a pilot set of relevant predictors; then we select the
bandwidth $\hat b_n$ by applying the GCV method to this pilot set.

\subsection{Locally stationary autoregressive processes}\label{subsecTVAR}
Modeling a nonstationary process by autoregressive models with
time-varying coefficients has attracted considerable attention. A
traditional approach is to project the coefficient function onto a
basis of temporal functions, and estimates the basis coefficients;
see, for example, \cite
{SubbaRao1970,Grenier1983,RajanRaynerGodsill1997}. Other
contributions on parameter estimation can be found
in \cite
{Dahlhaus1997,DahlhausNeumannSachs1999,MoulinesPriouretRoueff2005,GencagaKuruogluErtuzunYildirim2008}
among others. Abramovich et al.~\cite{AbramovichSpencerTurley2007}
considered the problem of order selection by requiring multiple
realizations. Let $a_1(\cdot), \ldots, a_p(\cdot)$ be
continuous\vadjust{\goodbreak}
functions. We shall prove that the time-varying autoregressive process
%
\begin{equation}
\label{eqnTVAR} y_i = a_1(i/n) y_{i-1} +
\cdots+ a_p(i/n) y_{i-p} + e_i,\qquad i = 1,\ldots,n,
\end{equation}
has an approximate solution of form (\ref{eqnformulation}),
and the difference is of a negligible order. Hence the results in
Section \ref{secmain} can be directly applied to address the
problem of parameter estimation, hypothesis testing and order
selection for time-varying autoregressive models.

Recall that $e_i = H(i/n; \bolds{\mathcal{F}}_i)$. Let
$\mathbf{x}_i = (y_{i},\ldots,y_{i-p+1})^\top$,
\[
\mathbf{A}(t) =
\pmatrix{ a_1(t)
& \cdots& a_{p-1}(t) & a_p (t)
\vspace*{2pt}\cr
1 & \cdots& 0 & 0
\vspace*{2pt}\cr
\vdots& \ddots& \vdots& \vdots
\vspace*{2pt}\cr
0 & \cdots& 1 & 0 }
 \in\mathbb
R^{p \times p}
\]
and
\[
\mathbf{H}^\diamond(t;\bolds{\mathcal F}_k) =
\pmatrix{
H(t;\bolds{\mathcal F}_k)
\vspace*{2pt}\cr
0
\vspace*{2pt}\cr
\vdots
\vspace*{2pt}\cr
0 } \in\mathbb R^p.
\]
Then (\ref{eqnTVAR}) can be written as
%
\begin{equation}
\label{eqnTVARvector} \mathbf{x}_i = \mathbf{A}(i/n)
\mathbf{x}_{i-1} + \mathbf{H}^\diamond(i/n;\bolds{\mathcal
F}_i).
\end{equation}
We shall make the following assumptions:
\begin{longlist}[(T1)]
\item[(T1)] The starting point $(y_p,\ldots,y_1)^\top\in
\mathcal{L}^2$.

\item[(T2)] The coefficient functions
$a_j(\cdot)$, $j = 1,\ldots,p$, are Lipschitz continuous on $[0,1]$.

\item[(T3)] $\sum_{j=1}^p a_j(t) z^j \neq1$ for all $|z| \leq
1+c$ with $c > 0$ uniformly in $t \in[0,1]$.
\end{longlist}
Conditions (T2) and (T3) entail local stationarity and short-range
dependence, respectively; see also \cite{Dahlhaus1996}. Proposition
\ref{propTVAR} states that the autoregressive process~(\ref{eqnTVARvector}) can be approximated by
(\ref{eqnformulation}) with a uniform approximation error of
order $O_p(n^{-1})$.

%
\begin{proposition}\label{propTVAR}
Assume \textup{(T1)--(T3)}. If $H \in \operatorname{Lip}$, then there exists a~measurable
function $\mathbf{G} \in \operatorname{Lip}$ and a constant $C > 0$ such that
%
\begin{equation}
\label{eqnTVARmaxdiff} \max_{1 \leq i \leq n} \bigl\|\mathbf{x}_i -
\mathbf{G}(i/n; \bolds{\mathcal{F}}_i)\bigr\| \leq C n^{-1}.
\end{equation}
\end{proposition}

In proving asymptotic results of Section \ref{secmain}, the key
quantity is the partial sum process $\sum_{i=1}^k \mathbf{x}_i
\mathbf{x}_i^\top$ and $\sum_{i=1}^k \mathbf{x}_i e_i$, $k =
1,\ldots,n$. By Proposition \ref{propTVAR}, there exists a
measurable function $\mathbf{G} \in \operatorname{Lip}$ such that
\[
\max_{1 \leq k \leq n} \Biggl\llvert \sum_{i=1}^k
\bigl\{\mathbf{x}_i \mathbf{x}_i^\top-
\mathbf{G}(i/n;\bolds{\mathcal F}_i) \mathbf{G}(i/n;\bolds{\mathcal
F}_i)^\top\bigr\}\Biggr\rrvert = O_p(1)
\]
and
\[
\max_{1 \leq k \leq n} \Biggl\llvert \sum_{i=1}^k
\bigl\{\mathbf{x}_i - \mathbf{G}(i/n;\bolds{\mathcal
F}_i)\bigr\}e_i\Biggr\rrvert = O_p(1).
\]
A careful check of the proofs of our results in Section
\ref{secmain} indicates that, due to the above relation, they are
still valid for the time-varying autoregressive process
(\ref{eqnTVAR}).

\subsection{A comparison with GLRT}
\label{subsecGLRT}
The generalized likelihood ratio test (GLRT, \cite
{FanZhangZhang2001}) is a popular method for nonparametric
hypothesis testing. It was used by Cai, Fan and Li~\cite
{CaiFanLi2000} for testing the coefficient constancy, and
generalized by Fan and Huang \cite{FanHuang2005} to semiparametric
models. Properties of GLRT have been extensively studied for i.i.d.
samples. However, its validity for dependent data is not guaranteed. We
shall here briefly review the GLRT and compare it with our method. For
the null hypothesis
\[
H_0\dvtx  \bolds{\beta}(\cdot) \equiv\bolds{\beta}
\]
for some vector $\bolds{\beta}\in\mathbb R^p$, the GLRT statistic
is defined as
\[
T_{\mathrm{GLR}} = \frac{n }{2} \log\frac{\operatorname{RSS}_0 }{\operatorname
{RSS}_1} = \frac{n }{2}
\log\frac{\sum_{i=1}^n \{y_i - \mathbf{x}_i^\top
\breve{\bolds{\beta}}\}^2 }{\sum_{i=1}^n \{y_i - \mathbf{x}_i^\top\tilde{\bolds{\beta}}(i/n)\}^2},
\]
where $\breve{\bolds{\beta}}$ is the least squares estimate. To
construct the null distribution of $T_{\mathrm{GLR}}$, we use the
conditional bootstrap as suggested by \cite
{CaiFanLi2000,FanHuang2005}. Let $\tilde\sigma^2 = n^{-1}
\operatorname{RSS}_1$ and $\{e_i^\diamond\}_{i \in\mathbb Z}$ be i.i.d.
$N(0,\tilde\sigma^2)$. We generate the bootstrap sample $y_i^\diamond
= \mathbf{x}_i^\top\breve{\bolds{\beta}} + e_i^\diamond$, $i
= 1,\ldots,n$, and compute the test statistic $T_{\mathrm
{GLR}}^\diamond$. We approximate the distribution of $T_{\mathrm
{GLR}}$ by that of $T_{\mathrm{GLR}}^\diamond$.

\begin{figure}
\centering
\begin{tabular}{@{}cc@{}}

\includegraphics{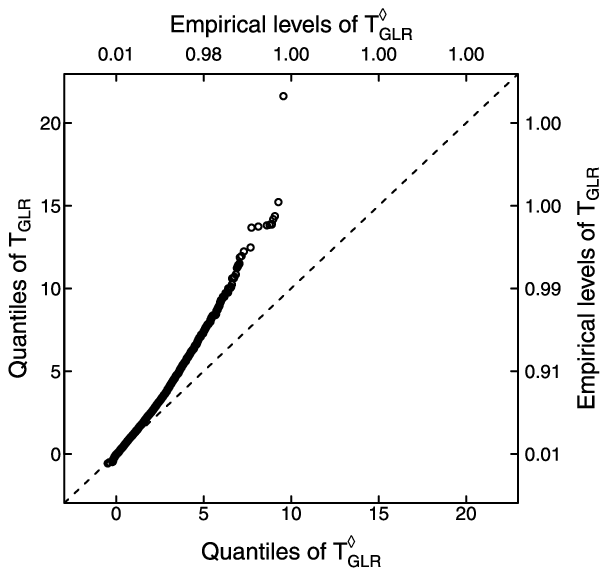}
 & \includegraphics{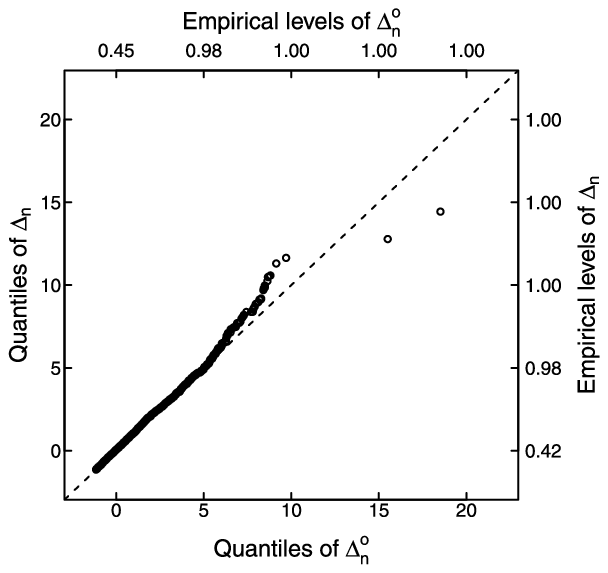}\\
\footnotesize{(a)} & \footnotesize{(b)}
\end{tabular}
\caption{A comparison of the GLRT \textup{(a)} with our dependence-adjusted
testing procedure~\textup{(b)}.
Q--Q plots of $T_{\mathrm{GLR}}$ against $T_{\mathrm
{GLR}}^\diamond$ \textup{(a)} and $\Delta_n(\mathbf{A},
\hat{\mathbf{a}},\mathbf{W})$ against $\Delta_n^\circ
(\mathbf{A}, \hat{\mathbf{a}},\mathbf{W})$~\textup{(b)}.
The dashed lines in \textup{(a)} and \textup{(b)} have unit slope and zero
intercept.}\label{figQQPlots}
\end{figure}

Consider the AR-ARCH process with time-varying conditional variance
\begin{eqnarray*}
y_i & = & 0.5 y_{i-1} + 0.25\bigl[ 1 + \bigl\{1+
\exp(3-6i/n)\bigr\}^{-1}\bigr] e_i,
\\
e_i & = & \bigl(1 + 0.25 e_{i-1}^2
\bigr)^{1/2} \varepsilon_i,
\end{eqnarray*}
where $\{\varepsilon_i\}_{i \in\mathbb Z}$ are i.i.d. $N(0,1)$. Let $n =
500$ and the bandwidth $b_n = n^{-1/5} = 0.289$. We consider testing
whether the coefficient of $y_{i-1}$ is a constant. For $\Delta_n(\mathbf{A}, \hat{\mathbf{a}},\mathbf{W})$, we use the
identity weights $\mathbf{W} = \mathbf I_{p \times p}$ and
obtain its cut-off value by the simulation-assisted procedure in
Section \ref{subsechypotestsimubased} with 5000 simulated $\Delta_n^\circ(\mathbf{A}, \hat{\mathbf{a}},\mathbf{W})$. For
$T_{\mathrm{GLR}}$, the cut-off value is obtained by 5000
bootstrapped~$T_{\mathrm{GLR}}^\diamond$. We generate 5000
realizations of the AR-ARCH process and use Q--Q plots to examine the
performance. The results are presented in Figure~\ref{figQQPlots}. It
shows that the GLRT fails to provide valid $p$-values in the presence
of dependence and nonstationarity. For example, the empirical
acceptance probabilities are 79\%, 86.4\% and 94.7\% for the 90\%, 95\%
and 99\% nominal levels, respectively. As shown in Figure~\ref
{figQQPlots}(b), our dependence-adjusted procedure provides a
satisfactory approximation of $\Delta_n(\mathbf{A}, \hat
{\mathbf{a}}, \mathbf{W})$. At 90\%, 95\% and 99\% nominal
levels, our empirical acceptance probabilities are 89.5\%, 95.0\% and
98.7\%, respectively.

%

%
\begin{table}
\caption{Percentages of under-fitted, correctly fitted and over-fitted
models selected by the variable selection information criterion
(\protect\ref{eqnVIC})
for $n = 500$. Medians of the selected bandwidths are $\hat b_n(\mathrm{i}) =
0.25$ and $\hat b_n(\mathrm{ii}) = 0.18$ for models \textup{(i)} and \textup{(ii)},
respectively, and $c = 2/3$}\label{tabsimulationVIC}
\begin{tabular*}{\textwidth}{@{\extracolsep{\fill}}lcd{3.2}c@{}}
\hline
&  \multicolumn{3}{c@{}}{\textbf{VIC}} \\[-6pt]
& \multicolumn{3}{c@{}}{\hrulefill} \\
\multicolumn{1}{@{}l}{$\bolds{b}$}  & \textbf{Under-fitted} & \multicolumn{1}{c}{\textbf{Correctly fitted}} & \multicolumn{1}{c@{}}{\textbf{Over-fitted}} \\
\hline
\multicolumn{4}{@{}l}{Model (i)} \\
0.1  & 0.0 & 100.0 & 0.0 \\
$c\hat b_n(\mathrm{i})$  & 0.0 & 100.0 & 0.0 \\
0.2  & 0.0 & 100.0 & 0.0 \\
$\hat b_n(\mathrm{i})$ & 0.0 & 100.0 & 0.0 \\
0.3  & 0.0 & 100.0 & 0.0 \\
$\hat b_n(\mathrm{i})/c$  & 0.0 & 100.0 & 0.0 \\
0.4  & 0.0 & 100.0 & 0.0 \\
0.5  & 0.0 & 100.0 & 0.0 \\
0.6  & 0.0 & 100.0 & 0.0 \\
0.7  & 0.0 & 100.0 & 0.0 \\
0.8  & 0.0 & 100.0 & 0.0 \\
0.9  & 0.0 & 100.0 & 0.0 \\[3pt]
\multicolumn{4}{@{}l}{Model (ii)} \\
0.1  & 0.0 & 99.8 & 0.2 \\
$c\hat b_n(\mathrm{ii})$  & 0.0 & 100.0 & 0.0 \\
$\hat b_n(\mathrm{ii})$  & 0.0 & 100.0 & 0.0 \\
0.2  & 0.0 & 100.0 & 0.0 \\
$\hat b_n(\mathrm{ii})/c$  & 0.0 & 100.0 & 0.0 \\
0.3  & 0.0 & 100.0 & 0.0 \\
0.4  & 0.0 & 100.0 & 0.0 \\
0.5  & 0.0 & 100.0 & 0.0 \\
0.6  & 0.0 & 100.0 & 0.0 \\
0.7 & 0.2 & 99.8 & 0.0 \\
0.8  & 1.1 & 98.9 & 0.0 \\
0.9  & 3.3 & 96.7 & 0.0 \\
\hline
\end{tabular*}
\end{table}

\subsection{Simulation studies}\label{subsecsimulation}
We shall here carry out a simulation study to examine the
finite-sample performance of our hypothesis testing procedure in
Section \ref{subsechypotest} and the information criterion for
variable selection in Section \ref{subsecvarsel}. Let $P_j(t)$ be
the $j$th order Legendre polynomial and $\mathbf P(t) \in
\mathbb R^{5 \times5}$ be the diagonal matrix with $j$th
diagonal component $P_j(2t-1)/4$. Let $\bolds{\varepsilon}_k =
(\varepsilon_{k,1},\ldots,\varepsilon_{k,6})^\top$, $k \in\mathbb Z$,
be i.i.d. Rademacher random variables and $\mathbf{M}^\diamond=
(0.2^{|i-j|})_{1 \leq i,j \leq5}$. Then $\bolds\xi_k =
\mathbf{M}^\diamond(\varepsilon_{k,1}, \ldots,
\varepsilon_{k,5})^\top$, $k \in\mathbb Z$, forms a sequence of
independent random vectors with correlated components. Let
$\mathbf{x}_i = \sum_{j=0}^\infty\mathbf P(i/n)^j
\bolds\xi_{i-j}$ and $e_i = \sum_{j=0}^\infty P_6(i/n)^j
\varepsilon_{i-j,6}$. Consider:
\begin{longlist}[(ii)]
\item[(i)] a linear model with heteroscedastic errors: for $i =
1,\ldots,n$,
\[
y_i = (2i/n-1)^2 + 2 x_{i,1} + 2\log(i/n+1)
x_{i,2} + 0.5\bigl(x_{i,2}^2 +
x_{i,3}^2\bigr)^{1/2} e_i;
\]
\item[(ii)] a linear model with autoregressive effects: for $i =
1,\ldots,n$,
\[
y_i = 0.4 \sin(2 \pi i/n) y_{i-1} + 0.3 x_{i,1}
+ 0.4(2i/n-1)^3 x_{i,2} + \exp(0.5i/n-2)
\varepsilon_{i,6}.
\]
\end{longlist}
Let $n = 500$ and the bandwidth $b = 0.1k$, $k = 1,\ldots,9$. We
use the information criterion (\ref{eqnVIC}) to estimate the set
of relevant predictors. For model (ii), the whole set of potential
predictors is taken to be $(\mathbf{x}_i)$ along with three
lags of the response variable. A realization is categorized as
under-fitting if we miss at least one relevant predictor, and
over-fitting if the selected set contains at least one irrelevant
predictor without under-fitting. The results are summarized in
Table \ref{tabsimulationVIC} based on 5000 realizations. Given
models~(i) and~(ii), we use the hypothesis testing procedure to
test whether $(x_{i,1})$ has time-invariant contributions. Three
types of weight matrices are used: the identity
weights $\mathbf{W}_1(t) = \mathbf I_{s \times s}$, the
normalizer weights $\mathbf{W}_2(t) = \{\mathbf{A}
\bolds\Xi(t) \mathbf{A}^\top\}^{-1}$ and the prediction
weights $\mathbf{W}_3(t) = \mathbf{A} \mathbf{M}(\mathbf{G},t) \mathbf{A}^\top$. For each configuration, we
use the simulation-assisted hypothesis testing procedure in
Section \ref{subsechypotestsimubased} to obtain cut-off values
$\hat{q}_{0.90}$ and $\hat{q}_{0.95}$ with 5000 simulated $\hat
\Delta_n^\circ(\mathbf{A}, \hat{\mathbf{a}},\mathbf{W})$.
We then generate 5000 realizations of both models (i) and (ii),
and calculate the corresponding test statistic $\hat
\Delta_n(\mathbf{A}, \hat{\mathbf{a}},\mathbf{W})$.
Empirical acceptance probabilities are reported in Table
\ref{tabsimulationTesting}.

%
\begin{table}
\caption{Empirical acceptance probabilities (in percentage) of the
simulation-assisted hypothesis testing procedure in
Section \protect\ref{subsechypotestsimubased} for $n = 500$. Medians
of the selected bandwidths are $\hat b_n(\mathrm{i}) = 0.25$
and $\hat b_n(\mathrm{ii}) = 0.18$ for models \textup{(i)} and \textup{(ii)}, respectively, and $c
= 2/3$}\label{tabsimulationTesting}
\begin{tabular*}{\textwidth}{@{\extracolsep{\fill}}lcccccc@{}}
\hline
& \multicolumn{2}{c}{$\mathbf{W_{1}}\bolds{(t)}$}  & \multicolumn{2}{c}{$\mathbf{W_2}\bolds{(t)}$}  & \multicolumn{2}{c}{$\mathbf{W_3}\bolds{(t)}$} \\[-6pt]
& \multicolumn{2}{c}{\hrulefill}  & \multicolumn{2}{c}{\hrulefill}  & \multicolumn{2}{c@{}}{\hrulefill} \\
$\bolds{b}$ & \textbf{90\%} & \textbf{95\%} & \textbf{90\%} & \textbf{95\%} &\textbf{90\%} & \textbf{95\%} \\
\hline
\multicolumn{7}{@{}l}{Model (i)} \\
0.1 &  92.7 & 96.9 &  93.2 & 97.1 &  92.4 & 96.5 \\
$c\hat b_n(\mathrm{i})$ &  92.2 & 96.5 &  92.9 & 96.7 &  92.0 & 96.4 \\
0.2 &  91.7 & 96.0 &  92.4 & 96.2 &  91.7 & 95.9 \\
$\hat b_n(\mathrm{i})$ &  90.9 & 95.8 &  91.4 & 96.4 &  90.3 & 95.5 \\
0.3 &  90.3 & 95.2 &  90.9 & 95.7 &  90.0 & 94.8 \\
$\hat b_n(\mathrm{i})/c$ &  90.5 & 95.1 &  90.9 & 95.8 &  90.3 & 95.1 \\
0.4 &  90.9 & 95.5 &  91.7 & 96.1 &  91.1 & 95.5 \\
0.5 &  90.4 & 94.8 &  91.3 & 95.2 &  90.4 & 95.0 \\
0.6 &  90.7 & 95.3 &  91.2 & 95.5 &  90.7 & 95.2 \\
0.7 &  91.6 & 95.6 &  91.7 & 95.7 &  91.5 & 95.7 \\
0.8 &  89.1 & 94.9 &  89.0 & 94.8 &  89.0 & 94.8 \\
0.9 &  89.9 & 95.0 &  90.4 & 95.4 &  89.8 & 94.9 \\[3pt]
\multicolumn{7}{@{}l}{Model (ii)} \\
0.1 &  93.1 & 97.0 &  93.2 & 97.0 &  92.9 & 96.9 \\
$c\hat b_n(\mathrm{ii})$ &  92.6 & 96.6 &  92.2 & 96.1 &  92.1 & 96.2 \\
$\hat b_n(\mathrm{ii})$ &  91.1 & 96.1 &  91.0 & 95.5 &  90.7 & 95.8 \\
0.2 &  91.3 & 96.3 &  91.2 & 96.2 &  91.0 & 96.0 \\
$\hat b_n(\mathrm{ii})/c$ &  90.9 & 95.4 &  90.8 & 95.7 &  90.7 & 95.2 \\
0.3 &  90.8 & 95.1 &  89.9 & 94.8 &  90.4 & 95.1 \\
0.4 &  90.0 & 95.1 &  89.8 & 94.9 &  89.9 & 95.1 \\
0.5 &  89.5 & 94.6 &  88.2 & 93.9 &  89.2 & 94.5 \\
0.6 &  88.6 & 93.6 &  87.7 & 93.1 &  88.5 & 93.6 \\
0.7 &  88.6 & 93.8 &  87.4 & 93.6 &  88.1 & 93.8 \\
0.8 &  87.4 & 92.6 &  87.0 & 92.2 &  87.4 & 92.7 \\
0.9 &  88.3 & 94.4 &  88.1 & 94.0 &  88.2 & 94.3 \\
\hline
\end{tabular*}
\end{table}

%
\begin{figure}[b]

\includegraphics{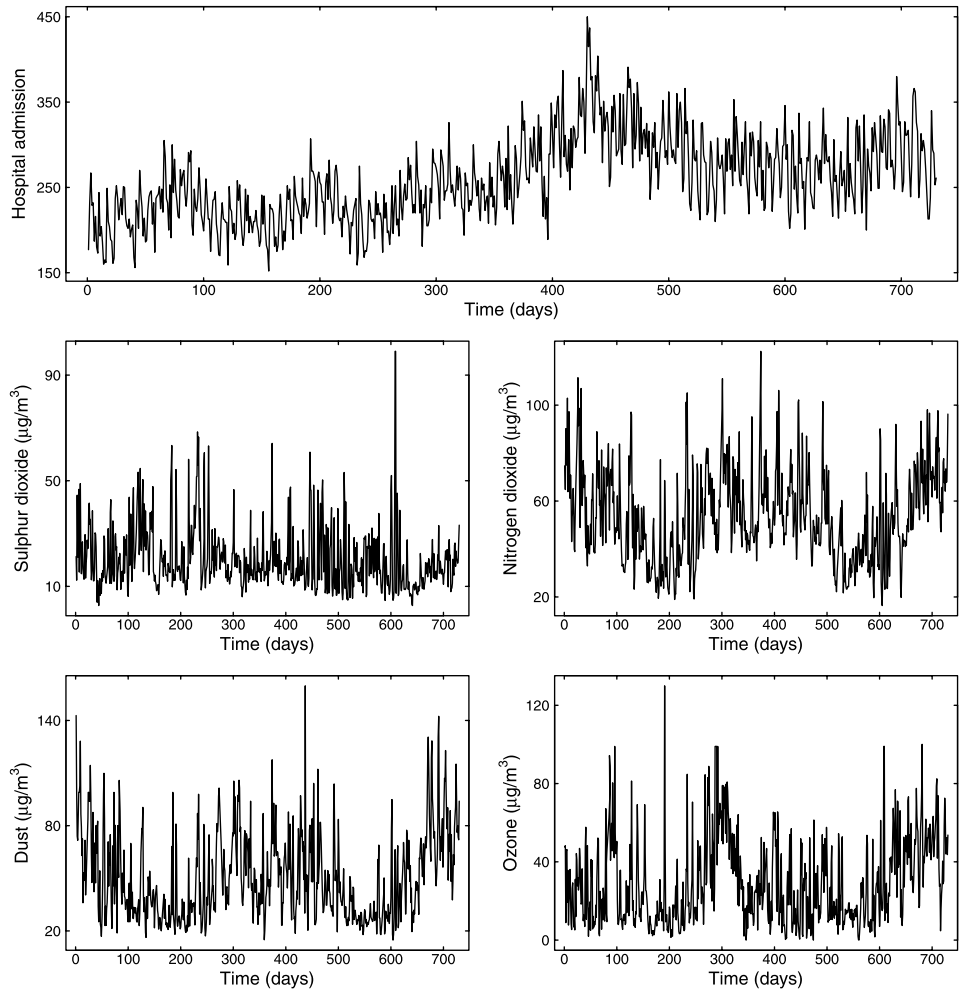}

\caption{Time series plots for daily hospital admission (top) and
levels of sulphur dioxide (middle left),
nitrogen dioxide (middle right), dust (bottom left) and ozone (bottom
right) from January 1, 1994 to December 31, 1995.}\label{figTSPlots}
\end{figure}

It can be seen that the empirical acceptance probabilities are
fairly close to their nominal levels (90\% and 95\%), and the
information criterion (\ref{eqnVIC}) performs quite well. In
addition, the results are not sensitive to choices of weight
matrices and bandwidths. For models (i) and (ii), medians of the
selected bandwidths based on the GCV criterion (\ref{eqnGCV}) are
$\hat b_n\mathrm{(i)} = 0.25$ and $\hat b_n(\mathrm{ii}) = 0.18$, respectively.
Observe that for model (i), the performance is also quite
satisfactory if we choose bandwidths $0.25 c$, 0.25 and $0.25/c$
with $c = 2/3$. A similar claim can be made for model (ii) as
well.\eject

\subsection{A real-data example}\label{subsecapplication}
\begin{figure}

\includegraphics{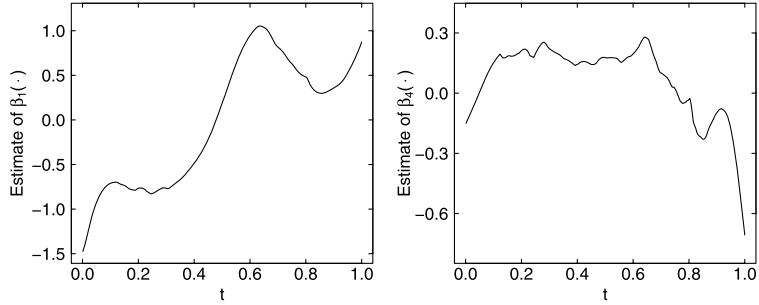}

\caption{Plots of estimated coefficient functions,
$\beta_1(\cdot)$ (left) and $\beta_4(\cdot)$ (right).}
\label{figbetatilde}
\end{figure}

We apply our model selection method to the Hong Kong circulatory and
respiratory data which contains
daily measurements of pollutants and hospital admissions in Hong
Kong between January~1, 1994 and December 31, 1995 ($n = 730$).
Four pollutants, sulphur dioxide (in $\upmu \mathrm{g}/\mathrm{m}^3$), nitrogen
dioxide (in $\upmu \mathrm{g}/\mathrm{m}^3$), dust (in $\upmu \mathrm{g}/\mathrm{m}^3$) and ozone (in $\upmu
\mathrm{g}/\mathrm{m}^3$), are considered here. The purpose is to understand the
association between daily hospital admission ($y_i$) and levels of
sulphur dioxide ($x_{i,2}$), nitrogen dioxide ($x_{i,3}$), dust
($x_{i,4}$) and ozone ($x_{i,5}$). Figure \ref{figTSPlots}
provides their time series plots. In the analysis, we regularize
the data so that each variable has zero mean and unit variance.
Letting $x_{i,1} \equiv1$ be the intercept, we consider the
time-varying coefficient model
%
\begin{equation}
\label{eqnHKData} y_i = \beta_1(i/n) + \sum
_{j = 2}^5 \beta_j(i/n)
x_{i,j} + e_i, \qquad i = 1,\ldots,n.
\end{equation}
The dataset has been studied by \cite
{FanZhang1999,FanZhang2000,CaiFanLi2000} by assuming that the
observations are i.i.d., while Zhou and Wu \cite{ZhouWu2010} found
substantial dependence among the fitted residuals. We shall here model
the process by
(\ref{eqnformulation}) and apply our model selection method in
Section \ref{secmain}. The selected bandwidth and tuning
parameter are $\hat b_n = 0.13$ and $\hat\chi_n = 0.072$,
respectively. The information criterion (\ref{eqnVIC}) selects
the intercept ($x_{i,1}$), nitrogen dioxide ($x_{i,3}$) and dust
($x_{i,4}$) as relevant predictors. Fan and Zhang~\cite
{FanZhang2000} did not
consider the ozone effect ($x_{i,5}$) and concluded that sulphur
dioxide ($x_{i,2}$) is not statistically significant. We then
apply the hypothesis testing procedure in Section
\ref{subsechypotestsimubased} to examine whether the selected
variables really have time-varying contributions. With 5000
simulated $\Delta_n^\circ(\mathbf{A}, \hat{\mathbf{a}},\mathbf{W})$, the results are summarized in Table
\ref{tabapplication}. Hence, at 10\% significance level, we
conclude that $\beta_1(\cdot)$ and $\beta_4(\cdot)$ are
time-varying while $\beta_3(\cdot)$ can be treated as
time-invariant, suggesting the model
%
\begin{equation}
\label{eqnHKDataselected} y_i = \beta_1(i/n) +
\beta_3 x_{i,3} + \beta_4(i/n)
x_{i,4} + e_i,\qquad i = 1,\ldots,n,
\end{equation}
where $\hat\beta_3 = \int_0^1 \tilde\beta_3(t) \,dt = 0.15$, and
$\tilde\beta_1(\cdot)$ and $\tilde\beta_4(\cdot)$ are plotted in
Figure \ref{figbetatilde}.

%
\begin{table}
\caption{Summary of test statistics and corresponding $p$-values for
testing parameter constancy with 5000 simulated $\Delta_n^\circ
(\mathbf{A}, \hat{\mathbf{a}},\mathbf{W})$}\label{tabapplication}
\begin{tabular*}{\textwidth}{@{\extracolsep{\fill}}ld{2.2}cd{3.2}cd{2.2}c@{}}
\hline
&  \multicolumn{2}{c}{$\mathbf{W_1}\bolds{(t)}$} &  \multicolumn{2}{c}{$\mathbf{W_2}\bolds{(t)}$} &
\multicolumn{2}{c}{$\mathbf{W_3}\bolds{(t)}$} \\[-6pt]
&  \multicolumn{2}{c}{\hrulefill} &  \multicolumn{2}{c}{\hrulefill} &
\multicolumn{2}{c@{}}{\hrulefill} \\
& \multicolumn{1}{c}{$\bolds{\Delta_n}$} & \multicolumn{1}{c}{$\bolds{p}$\textbf{-value}} & \multicolumn{1}{c}{$\bolds{\Delta_n}$} &
\multicolumn{1}{c}{$\bolds{p}$\textbf{-value}} & \multicolumn{1}{c}{$\bolds{\Delta_n}$}
& \multicolumn{1}{c@{}}{$\bolds{p}$\textbf{-value}} \\
\hline
$\beta_1(\cdot)$ & 69.77 & 0.00 & 120.77 & 0.02 & 69.77 & 0.00
\\
$\beta_3(\cdot)$ & 6.88 & 0.14 & 12.47 &
0.19 & 6.85 & 0.15 \\
$\beta_4(\cdot)$ & 16.27 & 0.02 & 30.13 & 0.09 & 23.06 & 0.01 \\
\hline
\end{tabular*}
\end{table}

%

\begin{appendix}

\section*{Appendix}\label{secproof}
For $a, b \in\mathbb R$, write $a \wedge b =
\min(a, b)$ and $a \vee b = \max(a, b)$. For a matrix $\mathbf{A}$, recall that $\underline\rho(\mathbf{A}) = \inf
\{|\mathbf{A} \mathbf{v}|\dvtx  |\mathbf{v}| = 1 \}$ and
$\overline\rho(\mathbf{A}) = \sup\{|\mathbf{A}
\mathbf{v}|\dvtx  |\mathbf{v}| = 1 \}$. The proofs of the
following two propositions are straightforward, and the details are
omitted.

%
\setcounter{proposition}{0}
\begin{proposition}
\label{propmatrixnorm} Let $\mathbf{A} = (a_{ij})_{1 \leq i
\leq I, 1 \leq j \leq J}$ be a real matrix.  Then\break \textup{(i)}~$\max_{i,j}
|a_{ij}| \le\overline{\rho} (\mathbf{A}) \le\sqrt{I J}
\max_{i,j} |a_{ij}|$; \textup{(ii)} If $\mathbf B$ has same dimension
as~$\mathbf{A}$, then $\overline{\rho}(\mathbf{A} +
\mathbf B) \leq\overline{\rho}(\mathbf{A}) +
\overline{\rho}(\mathbf B)$; \textup{(iii)} If $\mathbf B = (b_{j
k})_{1 \leq j \leq J, 1 \leq k\leq K} \in\mathbb{R}^{J \times
K}$, then $\overline{\rho}(\mathbf{A} \mathbf B) \leq
\overline{\rho}(\mathbf{A}) \overline{\rho}(\mathbf B)$ and
$\underline{\rho}(\mathbf{A} \mathbf B) \geq
\underline{\rho}(\mathbf{A}) \underline{\rho}(\mathbf B)$;
and \textup{(iv)} $\overline{\rho}(\mathbf{a} \mathbf{a}^\top) =
|\mathbf{a}|^2$ for any column vector~$\mathbf{a}$.
\end{proposition}

%
\begin{proposition}
\label{propinversenorm} Assume that $\mathbf{A}$ is a
nonsingular square matrix and that $\mathbf E$ is a matrix
with the same dimension. If $\overline\rho(\mathbf{A}^{-1}
\mathbf E) < 1$, then $\mathbf{A} + \mathbf E$ is
nonsingular and $\overline\rho\{(\mathbf{A} + \mathbf
E)^{-1} - \mathbf{A}^{-1}\} \leq\overline\rho(\mathbf E)
\overline\rho(\mathbf{A}^{-1})^2 / \{1 - \overline
\rho(\mathbf{A}^{-1} \mathbf E)\}$.
\end{proposition}

Let $\bolds{\mathcal F}_{i,j} = (\bolds{\varepsilon}_i,
\ldots, \bolds{\varepsilon}_j)$, $i \leq j$. Define the
projection operator
\[
\mathcal{P}_k \cdot= E(\cdot| \bolds{\mathcal{F}}_k) -
E(\cdot| \bolds{\mathcal{F}}_{k-1}),\qquad  k \in\mathbb{Z}.
\]
Let $\bolds\vartheta_k(t) = \mathbf{J}(t; \bolds
{\mathcal{F}}_k)$ be a zero mean process. Write $t_{i,n} = i/n$, $i =
1,\ldots,n$. Lemmas \ref{lembndlin} and \ref{lembndquad} provide
$\mathcal{L}^q$-bounds for linear and quadratic forms of $\{
\bolds\vartheta_k(t_{k,n})\}_{k=1}^n$, respectively. To prove
Theorem \ref{thmCLTahat}, we need Lemmas \ref{lembndlin} and \ref
{lemsupbnd}.

%
\begin{lemma}\label{lembndlin}
Assume $\Theta_{0,q}(\mathbf{J}) < \infty$, $q > 1$. Write $q'
= q \wedge2$. Let\break $\{\mathbf{A}_{k,n}(t)\}_{k = 1}^n$,\vspace*{1pt} $t \in
[0,1]$, be a
sequence of real matrix functions, and define $\mathbf S_n(t) =
\sum_{k=1}^n \mathbf{A}_{k,n}(t) \bolds
\vartheta_k(t_{k,n})$. Then:
\begin{longlist}[(ii)]
\item[(i)] $\|\mathbf S_n(t)\|_q \leq C_q
[\sum_{k=1}^n |\overline{\rho}\{\mathbf{A}_{k,n}(t)\}|^{q'}]^{1/q'} \Theta_{0,q}(\mathbf{J})$;
\item[(ii)] $\|\sup_{t \in[0,1]} |\mathbf S_n(t)|\|_q \leq C_q n^{1/q'}
\mathcal{A}_n \Theta_{0,q}(\mathbf{J})$,
\end{longlist}
where $\mathcal{A}_n =
\sup_{t \in[0,1]} [\overline{\rho}\{\mathbf{A}_{1,n}(t)\} +
\sum_{k=1}^{n-1} \overline{\rho}\{\mathbf{A}_{k+1,n}(t) -
\mathbf{A}_{k,n}(t)\}]$.
\end{lemma}

\begin{pf}
Let $\mathbf D_{k,l,n} = E
\{\bolds{\vartheta}_{k}(t_{k,n})
|\bolds{\mathcal{F}}_{k-l,k}\} - E
\{\bolds{\vartheta}_{k}(t_{k,n})
|\bolds{\mathcal{F}}_{k-l+1,k}\}$. Then
$\bolds{\vartheta}_k(t_{k,n}) - E
\bolds{\vartheta}_k(t_{k,n}) = \sum_{l=0}^{\infty}
\mathbf D_{k,l,n}$ and $\mathbf D_{k,l,n}$, $k=1,\ldots,
n$, form martingale differences. By the Burkholder and the
Minkowski inequalities, we have
\[
\Biggl\llVert \sum_{k=1}^n
\mathbf{A}_{k,n}(t) \mathbf D_{k,l,n}\Biggr\rrVert_q^{q'}
\leq C_q \sum_{k=1}^n \bigl|
\overline{\rho}\bigl\{\mathbf{A}_{k,n}(t)\bigr\}\bigr|^{q'} \|
\mathbf D_{k,l,n}\|_q^{q'}.
\]
Since $\|\mathbf D_{k,l,n}\|_q = \|E
\{\bolds{\vartheta}_{l}(t_{k,n})|\bolds{\mathcal
{F}}_{0,l}\}
- E \{\bolds{\vartheta}_{l,\{0\}}(t_{k,n})
|\bolds{\mathcal{F}}_{0,l}\} \|_q \leq
\delta_{l,q}(\mathbf{J})$, (i) follows. We now prove (ii).\vadjust{\goodbreak} By
Doob's inequality and the summation by parts formula, we have $\|
\sup_{t \in[0,1]} |\sum_{k=1}^n \mathbf{A}_{k,n}(t)
\mathbf D_{k,l,n}| \|_q \leq C_q \mathcal{A}_n n^{1/q'}
\delta_{l,q}(\mathbf{J})$, entailing~(ii).
\end{pf}

%
\begin{lemma}\label{lembndquad}
Assume $\Theta_{0,2q}(\mathbf{J}) < \infty$, $q \geq2$. Let
$\{\mathbf Q_{i,j,n}\}_{1 \leq i < j \leq n}$ be real matrices
and $L_n = \sum_{1 \leq i < j \leq n} \bolds
\vartheta_i(t_{i,n})^\top\mathbf Q_{i,j,n} \bolds
\vartheta_j(t_{j,n})$. Then
\[
\bigl\| L_n - E (L_n) \bigr\|_q \leq C_q
n^{1/2} \mathcal{Q}_n \Theta_{0,2q}(
\mathbf{J})^2,
\]
where $\mathcal{Q}_n^2 = (\max_i \sum_{j=i+1}^n
|\overline{\rho}(\mathbf Q_{i,j,n})|^2) \vee(\max_j
\sum_{i=1}^{j-1} |\overline{\rho}(\mathbf Q_{i,j,n})|^2)$.
\end{lemma}

\begin{pf}
Let $\tilde{\bolds{\vartheta}}_k(t) = E
\{\bolds{\vartheta}_k(t) |
\bolds{\mathcal{F}}_{k-m,k}\}$ be the $m$-dependent
approximated process and $\tilde L_n$ be the corresponding
quadratic form. If $l > 2m$, $\mathcal{P}_{j-l}
\{\tilde{\bolds{\vartheta}}_i(t_{i,n})^\top\times  Q_{i,j,n}
\tilde{\bolds{\vartheta}}_j(t_{j,n})\} = 0$. Hence
\[
\bigl\|\tilde L_n - E (\tilde L_n)\bigr\|_q \leq\sum
_{l=0}^{2m} \Biggl\llVert \sum
_{j=2}^n \mathcal{P}_{j-l} \sum
_{i=1}^{j-1} \tilde{\bolds{\vartheta}}_i(t_{i,n})^\top
Q_{i,j,n} \tilde{\bolds{\vartheta}}_j(t_{j,n})
\Biggr\rrVert_q,
\]
where
\[
\Biggl\llVert \sum_{j=2}^n
\mathcal{P}_{j-l} \Biggl(\sum_{i = 1}^{j-l-1}
+ \sum_{i = j-l}^{j-1} \Biggr) \tilde{\bolds{
\vartheta}}_i(t_{i,n})^\top Q_{i,j,n}
\tilde{\bolds{\vartheta}}_j(t_{j,n}) \Biggr
\rrVert_q^2 \leq C_q n \mathcal{Q}_n^2
l^2 \Theta_{0,2q}(\mathbf{J})^2.
\]
By Lemma \ref{lembndlin} and the arguments of Proposition 1 in \cite
{LiuWu2010}, we have $\|\{L_n - E (L_n)\} - \{\tilde L_n - E(\tilde
L_n)\}\|_q \leq C_q \sqrt{n} \mathcal{Q}_n \Theta_{0,2q}(\mathbf{J})^2$. So Lemma \ref{lembndquad} follows.
\end{pf}

%
\begin{lemma}\label{lemsupbnd}
Assume $\sup_{t \in[0,1]} \|\mathbf{J}(t;\bolds{\mathcal
F}_0)\|_\iota< \infty$, $\iota> 2$, and $\Theta_{n,
\iota}(\mathbf{J}) = O(n^{-\nu})$ for some $\nu
> 1/2 - 1/\iota$. Let $\mathbf S_{K,n}(t) = (n b_n)^{-1}
\sum_{k=1}^n K\{(k/n-t)/b_n\} \times \bolds\vartheta_k(t_{k,n})$.
Then
%
\setcounter{equation}{0}
\begin{equation}
\label{eqO19859a} \sup_{b_n \le t \le1-b_n} \bigl|\mathbf S_{K,n}(t)\bigl| =
\frac{ {O_p(v_n)}
}{{ n b_n}} \qquad\mbox{where } v_n = n^{1/\iota} + (n
b_n \log n)^{1/2}.
\end{equation}
\end{lemma}

\begin{pf}
Let $S_n^* = n b_n \sup_{b_n \le t \le1-b_n} |\mathbf
S_{K,n}(t)|$. By Theorem 2(ii) in \cite{LiuXiaoWu2012}, there
exist constants $C_1, C_2 > 0$ such that, for all $\lambda\ge1$ and $l$,
%
\begin{eqnarray}
\label{eqnmaxinq}
&&\operatorname{pr} \Biggl(\max_{0 \leq j \le n b_n} \Biggl\llvert \sum
_{i=l}^{l+j} \bolds\vartheta_i(t_{i,n})
\Biggr\rrvert \geq\lambda v_n \Biggr)
\nonumber
\\[-8pt]
\\[-8pt]
\nonumber
&&\qquad\leq C_1
\frac{n b_n }{(\lambda v_n)^\iota} + 2 \exp\bigl\{-(\lambda v_n)^2/(n
b_n C_2)\bigr\}.
\end{eqnarray}
Note that $[b_n, 1- b_n] \subseteq\bigcup_{j \leq1/b_n} [jb_n,
(j+1)b_n]$. Using the summation by parts formula, since $K$ has
support $[-1, 1]$, we have (\ref{eqO19859a}) in view of
(\ref{eqnmaxinq}) and
\[
\operatorname{pr}\bigl(S_n^* \ge\lambda v_n\bigr) = O
\bigl(b_n^{-1}\bigr) \biggl[ C_1
\frac{n b_n }{(\lambda v_n)^\iota} + 2 \exp\bigl\{-(\lambda v_n)^2/(n
b_n C_2)\bigr\} \biggr]
\]
by choosing a sufficiently large $\lambda$.\vadjust{\goodbreak}
\end{pf}

For $l \in\{0,1,2\}$, let
\[
\mathbf R_{n,l}(t) = (nb_n)^{-1} \sum
_{i=1}^n \mathbf{x}_i e_i
\bigl\{(i/n - t)/b_n\bigr\}^l K\bigl\{(i/n-t)/b_n
\bigr\}.
\]

\begin{pf*}{Proof of Theorem \protect\ref{thmCLTahat}}
By Lemma \ref{lembndlin}, we have
\begin{eqnarray*}
& & \int_0^1 \mathbf{A} \mathbf{M}(
\mathbf{G},t)^{-1} \mathbf R_{n,0}(t) \,dt
\\
& &\qquad=  \frac{1 }{ nb_n} \sum_{i=1}^n
\mathbf{A} \biggl\{\int_0^1 \mathbf{M}(
\mathbf{G},t)^{-1} K \biggl(\frac{i/n-t }{ b_n} \biggr) \,dt \biggr\}
\mathbf{x}_i e_i
\\
&&\qquad =  \frac{1 }{ n} \sum_{i=1}^n
\mathbf{A} \mathbf{M}(\mathbf{G},i/n)^{-1} \mathbf{x}_i
e_i + O_p \biggl\{\frac{(nb_n)^{1/2} }{ n} +
\frac{b_n n^{1/2} }{ n} \biggr\}.
\end{eqnarray*}
By $m$-dependence approximation, under Conditions (A2), (A3) and
(A4), we obtain
\[
n^{-1/2} \sum_{i=1}^n \mathbf{A}
\mathbf{M}(\mathbf{G},i/n)^{-1} \mathbf{x}_i
e_i \Rightarrow N \biggl\{0, \int_0^1
\mathbf{A} \bolds\Xi(t) \mathbf{A}^\top \,dt \biggr\}.
\]
By Lemmas \ref{lembndlin} and \ref{lemsupbnd}, and the argument
in the proof of Theorem 3 in~\cite{ZhouWu2010}, we have $\sup_{t
\in[0,1]} |\tilde{\bolds{\beta}}(t) - \bolds{\beta}(t)| =
O_p(\varphi_n)$ and
%
\begin{eqnarray}
\label{eqnBahadurrep}
&&\sup_{b_n \leq t \leq1-b_n} \bigl|\mathbf{M}(\mathbf{G},t)\bigl\{
\tilde{\bolds{\beta}}(t) - \bolds{\beta}(t) - 2^{-1}
b_n^2 \kappa_2 \bolds{\beta}''(t)
\bigr\} - \mathbf R_{n,0}(t)\bigr|
\nonumber
\\[-8pt]
\\[-8pt]
\nonumber
&&\qquad = O_p(\varphi_n
\rho_n).
\end{eqnarray}
Therefore,
\[
\hat{\mathbf{a}} - \mathbf{a} - \bolds\xi_n = \int
_0^1 \mathbf{A} \mathbf{M}(\mathbf{G},t)^{-1}
\mathbf R_{n,0}(t) \,dt + O_p\bigl(\varphi_n
\rho_n + b_n \varphi_n + b_n^3
\bigr).
\]
Under our bandwidth conditions, $\varphi_n \rho_n + b_n \varphi_n
+ b_n^3 = o(n^{-1/2})$. So Theorem~\ref{thmCLTahat} follows.
\end{pf*}

Let $\gamma_{k,2}(\mathbf{J}) = \sum_{i=0}^\infty
\delta_{i,2}(\mathbf{J}) \delta_{i+|k|,2}(\mathbf{J})$.
Lemma \ref{lemlrcov} provides continuity properties of long-run
covariance matrices for stochastically Lipschitz continuous
processes.

%
\begin{lemma}\label{lemlrcov}
Assume $\mathbf{J} \in \operatorname{Lip}$ and $\Theta_{0,2}(\mathbf{J}) <
\infty$. Then: \textup{(i)} for any nonnegative sequence $a_n \to0$, $\sup_{|t_1 - t_2| \leq a_n} \overline\rho\{\bolds{\Lambda}
(\mathbf{J}, t_1) - \bolds{\Lambda}(\mathbf{J}, t_2)\} =
o(1)$; \textup{(ii)} if, in addition, $\Theta_{n,2}(\mathbf{J}) = O(n^{-\nu
})$ for some $\nu> 0$, then $\sup_{|t_1 - t_2| \leq a_n} \overline
\rho\{\bolds{\Lambda}(\mathbf{J}, t_1) - \bolds{\Lambda}
(\mathbf{J}, t_2)\} = O\{a_n^{\nu/(1+\nu)}\}$; and \textup{(iii)} if $\inf_{t \in[0,1]} \underline\rho\{\bolds{\Lambda}(\mathbf{J},
t)\} > 0$, then \textup{(i)} and \textup{(ii)} hold for the inverse $\Lambda^{-1}(\mathbf{J}, t)$.
\end{lemma}

\begin{pf}
We first observe that
\begin{eqnarray*}
\overline\rho\bigl[E \bigl\{\bolds\vartheta_i(t_1)
\bolds\vartheta_j(t_2)^\top\bigr\}\bigr] &
\leq& \sum_{s \in\mathbb{Z}} \bigl\|\mathcal{P}_s \bolds
\vartheta_i(t_1) \bigr\| \bigl\|\mathcal{P}_s \bolds
\vartheta_j(t_2) \bigr\|
\\
& \leq& \sum_{s \in\mathbb{Z}} \delta_{i-s,2}(
\mathbf{J}) \delta_{j-s,2}(\mathbf{J}) = \gamma_{|j-i|,2}(
\mathbf{J}).
\end{eqnarray*}
The Lipschitz continuity implies
%
\begin{equation}
\label{eqnbndt1t2} \overline\rho\bigl(E \bigl[\bolds\vartheta_i(t_1)
\bigl\{\bolds \vartheta_j(t_2) - \bolds
\vartheta_j(t_1)\bigr\}^\top\bigr]\bigr) \leq C
\bigl\{\gamma_{|i-j|,2}(\mathbf{J}) \wedge|t_2 -
t_1|\bigr\}
\end{equation}
uniformly. Hence
\[
\sup_{|t_1 - t_2| \leq a_n} \overline\rho\bigl\{\bolds \Lambda(\mathbf{J},
t_1) - \bolds{\Lambda}(\mathbf{J}, t_2)\bigr\} \leq C
\sum_{k \in\mathbb{Z}} \bigl\{\gamma_{k,2}(\mathbf{J})
\wedge a_n\bigr\},
\]
which entails (i) by the dominated convergence theorem. Let $r_n =
a_n^{-1/(1+\nu)}$ which goes to infinity as $n \to\infty$. Since
$\sum_{k=l}^\infty
\gamma_{k,2}(\mathbf{J}) \leq\Theta_{l,2}(\mathbf{J})
\Theta_{0,2}(\mathbf{J})$, we have $\sum_{k=0}^\infty
\{\gamma_{k,2}(\mathbf{J}) \wedge a_n\} = \mathcal{O}(r_n a_n +
r_n^{-\nu})$, (ii) follows. Then (iii) follows by Proposition
\ref{propinversenorm}.
\end{pf}

Let $\mathbf{W}_0(\cdot)$ be a continuous mapping from $[0,1]$
to symmetric matrices in $\mathbb{R}^{p \times p}$. For $l \in
\{1,2\}$, define $\Lambda_{\mathbf{W}_0,l} =
\operatorname{tr}[\int_0^1 \{\mathbf{W}_0(t) \bolds
\Lambda(\mathbf{L},t)\}^l \,dt]$. Before we prove Theorem~\ref{thmCLTTn}, we shall first establish a parallel result for
%
\begin{equation}
\label{eqnTndiamond} T_n^\diamond(\mathbf{W}_0) =
\int_0^1 \mathbf R_{n,0}(t)^\top
\mathbf{W}_0(t) \mathbf R_{n,0}(t) \,dt.
\end{equation}
Let $r_{1,n} = (nb_n)^{-1} \sum_{k=0}^\infty\{\gamma_{k,2}
(\mathbf{L}) \wedge b_n\}$ and $r_{2,n} = (nb_n)^{-1}
\sum_{k=0}^\infty[\{k/(nb_n)\} \wedge1] \gamma_{k,2}
(\mathbf{L})$.

%
\begin{lemma}\label{lemCLTTndiamond}
Assume $\mathbf{L} \in \operatorname{Lip}$ and $\Theta_{0,4}(\mathbf{L}) <
\infty$. If $b_n \to0$ and $nb_n^{3/2} \to\infty$, then
%
\begin{equation}
\label{eqnCLTTndiamond} nb_n^{1/2} \bigl[T_n^\diamond(
\mathbf{W}_0) - E \bigl\{T_n^\diamond(
\mathbf{W}_0)\bigr\}\bigr] \Rightarrow N \bigl(0, 4 K^*_2
\Lambda_{\mathbf{W}_0,2}\bigr),
\end{equation}
and
%
\begin{equation}
\label{eqnTndiamondMean}\qquad E \bigl\{T_n^\diamond(
\mathbf{W}_0)\bigr\} = (nb_n)^{-1} K^*(0)
\Lambda_{\mathbf{W}_0,1} + O(r_{1,n} + r_{2,n}) + o
\bigl(n^{-1}b_n^{-1/2}\bigr).
\end{equation}
\end{lemma}

\begin{pf}
Let $\bolds\zeta_k(t) = \mathbf{L}(t;
\bolds{\mathcal{F}}_k)$ and $\tilde{\bolds{\zeta}}_k(t)
= E \{\bolds{\zeta}_k(t) |
\bolds{\mathcal{F}}_{k-m,k}\}$ be its $m$-dependent
counterpart and $\bolds{\Lambda}(\tilde{\mathbf{L}}, t)$ be
the corresponding long-run covariance matrix. Then $\bolds
\Lambda(\tilde{\mathbf{L}}, t) \to\bolds
\Lambda(\mathbf{L}, t)$ uniformly as $m \to\infty$. Let
$w_{k,n}(t) = (nb_n)^{-1}\times K\{(k/n - t)/b_n\}$, $k = 1,\ldots,n$,
and $\mathbf Q_{i,j,n} = \int_0^1 w_{i,n}(t) \mathbf{W}_0(t) w_{j,n}(t) \,dt$. The central limit theorem (\ref
{eqnCLTTndiamond}) is a
multivariate generalization of Theorem~A1 in \cite{ZhangWu2011} by
using Propositions \ref{propmatrixnorm} and
\ref{propinversenorm}. We shall only detail steps that require
special attention on the dimensionality. Essentially, we need to
show that
\begin{eqnarray*}
& & \lim_{m \to\infty} \lim_{n \to\infty} \Biggl[n^2b_n
\sum_{j=2m+1}^n \sum
_{i=1}^{j-2m} E\bigl\{\tilde{\mathbf
D}_{i,n}^{*\top} \mathbf Q_{i,j,n} E\bigl(\tilde{\mathbf
D}^*_{j,n} \tilde{\mathbf D}_{j,n}^{*\top}\bigr)
\mathbf Q_{i,j,n}^\top\tilde{\mathbf D}^*_{i,n}\bigr\}
\Biggr]
\\
&&\qquad =  K_2^* \Lambda_{\mathbf{W}_0,2},
\end{eqnarray*}
where $\tilde{\mathbf D}_{k,n}^* = \mathcal{P}_k
\sum_{l=0}^\infty\tilde{\bolds\zeta}_{k+l}(t_{k,n})$. Since
$E (\tilde{\mathbf D}_{k,n}^* \tilde{\mathbf
D}_{k,n}^{*\top}) = \bolds{\Lambda}(\tilde{\mathbf{L}},
t_{k,n})$, by Lem\-ma~\ref{lemlrcov}, we have
\begin{eqnarray*}
& & \sum_{j=2m+1}^n \sum
_{i=1}^{j-2m} E\bigl\{\tilde{\mathbf
D}_{i,n}^{*\top} \mathbf Q_{i,j,n} E\bigl(\tilde{\mathbf
D}^*_{j,n} \tilde{\mathbf D}_{j,n}^{*\top}\bigr)
\mathbf Q_{i,j,n}^\top\tilde{\mathbf D}^*_{i,n}\bigr\}
\\
&&\qquad =  \sum_{1 \leq i < j \leq n} \operatorname{tr}\bigl[\bigl\{
\mathbf{W}_0(t_{i,n}) \bolds{\Lambda}(\mathbf{L},
t_{i,n})\bigr\}^2\bigr] \biggl(\int_0^1
w_{i,n}(t) w_{j,n}(t) \,dt \biggr)^2
\\
& &\qquad\quad{} + o\bigl\{n^2b_n \bigl(n^2b_n
\bigr)^{-2}\bigr\} + O\bigl\{\ell(m) n^2b_n
\bigl(n^2b_n\bigr)^{-2}\bigr\} + O\bigl\{mn
\bigl(n^2b_n\bigr)^{-2}\bigr\}
\end{eqnarray*}
for some function $\ell(m) \to0$ as $m \to\infty$. Then
(\ref{eqnCLTTndiamond}) follows. For (\ref{eqnTndiamondMean}),
by the proof of Theorem 1 in \cite{ZhangWu2011}, we have
\[
E \bigl\{T_n^\diamond(\mathbf{W}_0)\bigr\} = \sum
_{i=1}^n \operatorname{tr}\bigl\{\mathbf
Q_{i,i,n} \bolds{\Lambda}(\mathbf{L}, t_{i,n})\bigr\} +
O(r_{1,n} + r_{2,n} + r_{3,n}),
\]
where
\[
r_{3,n} = \sum_{i=1}^n \Biggl(
\sum_{j = -\infty}^0 + \sum
_{j=n+1}^\infty \Biggr) \overline{\rho}(Q_{i,j,n})
\gamma_{|j-i|,2}(\mathbf{J}) \leq C nb_n \bigl(n^2b_n
\bigr)^{-1} \Theta_{0,2}(\mathbf{L})^2.
\]
Since $\sum_{i=1}^n \operatorname{tr}\{\mathbf Q_{i,i,n} \bolds
\Lambda(\mathbf{L}, t_{i,n})\} = (nb_n)^{-1} K^*(0)
\Lambda_{\mathbf{W}_0,1} + o(n^{-1} b_n^{-1/2})$,
(\ref{eqnTndiamondMean}) follows.
\end{pf}

\begin{pf*}{Proof of Theorem \ref{thmCLTTn}}
Let $\mathcal B_n = [0,b_n] \cup[1-b_n,1]$. Lemma
\ref{lembndlin} implies $\sup_{t \in\mathcal B_n} |\mathbf
R_{n,0}(t)| = O_p\{(nb_n)^{-1/2}\}$ and $\sup_{t \in\mathcal B_n}
|\mathbf{U}_n(t) - E\{\mathbf{U}_n(t)\}| =
O_p\{(nb_n)^{-1/2}\}$. By the proof of Theorem 1 in \cite
{ZhouWu2010}, we have $\sup_{t \in\mathcal B_n} |\tilde
{\bolds
\beta}(t) - \bolds{\beta}(t)| = O_p\{(nb_n)^{-1/2} + b_n^2\}$.
Hence,
\[
\int_{\mathcal B_n} \bigl\{\mathbf{A} \tilde{\bolds{\beta}}(t) -
\mathbf{a}\bigr\}^\top\mathbf{W}(t) \bigl\{\mathbf{A} \tilde{\bolds{
\beta}}(t) - \mathbf{a}\bigr\} \,dt = O_p\bigl(n^{-1} +
b_n^5\bigr).
\]
By (\ref{eqnBahadurrep}) and Lemma \ref{lemsupbnd}, we have
$\sup_{b_n \leq t \leq1-b_n}|\mathbf{A} \tilde{\bolds
\beta}(t) - \mathbf{a}| = O_p(\varphi_n)$ and
%
\begin{equation}
\label{eqnllestapp} \sup_{b_n \leq t \leq1-b_n} \bigl|\mathbf{A} \tilde{\bolds \beta}(t) -
\mathbf{a} - \mathbf{A} \mathbf{M}(\mathbf{G},t)^{-1} \mathbf
R_{n,0}(t)\bigr| = O_p(\varphi_n \rho_n).
\end{equation}
Let $\mathbf{W}_0(t) = \mathbf{M}(\mathbf{G},t)^{-1}
\mathbf{A}^\top\mathbf{W}(t) \mathbf{A} \mathbf{M}(\mathbf{G},t)^{-1}$. Since $\sup_{t \in\mathcal B_n}
|\mathbf R_{n,0}(t)| = \break O_p\{(nb_n)^{-1/2}\}$,
\[
\int_{\mathcal B_n} \mathbf R_{n,0}(t)^\top
\mathbf{W}_0(t) \mathbf R_{n,0}(t) \,dt = O_p
\bigl(n^{-1}\bigr).
\]
Under our bandwidth conditions, $nb_n^{1/2} \varphi_n^2 \rho_n =
o(1)$. For (\ref{eqnCLTTn}), by Lemma~\ref{lemCLTTndiamond}, it
suffices to show that both $r_{1,n}$ and $r_{2,n}$ are of order
$o(n^{-1}b_n^{-1/2})$. By the proof of Lemma \ref{lemlrcov},
$nb_n r_{1,n} = O\{b_n^{\nu/(1+\nu)}\} = o(b_n^{1/2})$ since $\nu
> 1$. Let $r_n = (nb_n)^{1/(2+\nu)}$. Then
\[
\sum_{k=0}^\infty\bigl[\bigl
\{k/(nb_n)\bigr\} \wedge1\bigr] \gamma_{k,2}(\mathbf{L})
\leq\frac{r_n(r_n+1) }{2nb_n} \Theta_{0,2}(\mathbf{L})^2 +
\Theta_{r_n,2}(\mathbf{L}) \Theta_{0,2}(\mathbf{L}) = O
\bigl(r_n^{-\nu}\bigr).
\]
Hence, we have $nb_n r_{2,n} = O\{(nb_n)^{-\nu/(2+\nu)}\} =
o(b_n^{1/2})$ since $\nu> 1$, (\ref{eqnCLTTn}) follows. Note
that
\[
T_n(\mathbf{A},\hat{\mathbf{a}},\mathbf{W}) - T_n(
\mathbf{A},\mathbf{a},\mathbf{W}) = I_n - 2 \II_n,
\]
where $I_n = \int_0^1 (\hat{\mathbf{a}} - \mathbf{a})^\top
\mathbf{W}(t) (\hat{\mathbf{a}} - \mathbf{a}) \,dt =
O_p(n^{-1})$, and by (\ref{eqnllestapp}) and Lemmas~\ref{lembndlin} and \ref{lemsupbnd},
\begin{eqnarray*}
\II_n & = & \int_0^1 (\hat{
\mathbf{a}} - \mathbf{a})^\top \mathbf{W}(t) \bigl\{\mathbf{A} \tilde{
\bolds{\beta}}(t) - \mathbf{a}\bigr\} \,dt
\\
& = & (\hat{\mathbf{a}} - \mathbf{a})^\top \biggl\{\int
_{b_n}^{1-b_n} \mathbf{W}(t) \mathbf{A} \mathbf{M}(
\mathbf{G},t)^{-1} \mathbf R_{n,0}(t) \,dt + O_p(
\varphi_n \rho_n + b_n \varphi_n)
\biggr\}
\\
& = & O_p\bigl\{(n^{-1/2}\bigl(n^{-1/2} +
\varphi_n \rho_n\bigr)\bigr\}.
\end{eqnarray*}
Note that $(nb_n)^{1/2} \varphi_n \rho_n = o(1)$, and Theorem
\ref{thmCLTTn} follows.
\end{pf*}

\begin{pf*}{Proof of Proposition \ref{thmPower}}
Under the local alternative (\ref{eqnlocalalternative}), we have
$\mathbf{A} \bolds{\beta}''(t) = d_n \mathbf{f}''(t)$
and
\[
T_n(\mathbf{A},\mathbf{a},\mathbf{W}) - T_n^\diamond(
\mathbf{W}_0) = d_n^2 \int
_0^1 \mathbf{f}(t)^\top\mathbf{W}(t)
\mathbf{f}(t) \,dt + I_n + 2\II_n,
\]
where by (\ref{eqnBahadurrep}) and Lemmas \ref{lembndlin} and
\ref{lemsupbnd},
\begin{eqnarray*}
I_n & = & \int_0^1 \bigl\{
\mathbf{A} \tilde{\bolds{\beta}}(t) - \mathbf{A} \bolds{\beta}(t)\bigr
\}^\top\mathbf{W}(t) \bigl\{ \mathbf{A} \tilde{\bolds{\beta}}(t) -
\mathbf{A} \bolds{\beta}(t)\bigr\} \,dt - T_n^\diamond(
\mathbf{W}_0)
\\
& = & O_p\bigl\{\bigl(\varphi_n + d_nb_n^2
\bigr) \bigl(\varphi_n\rho_n + d_nb_n^2
\bigr) + \bigl(\varphi_n\rho_n + d_nb_n^2
\bigr)\varphi_n\bigr\},
\end{eqnarray*}
the weight matrix $\mathbf{W}_0(t) = \mathbf{M}(\mathbf{G},t)^{-1} \mathbf{A}^\top\mathbf{W}(t) \mathbf{A}
\mathbf{M}(\mathbf{G},t)^{-1}$ and
\[
\II_n = d_n \int_0^1
\mathbf{f}(t)^\top\mathbf{W}(t) \bigl\{\mathbf{A} \tilde{\bolds{
\beta}}(t) - \mathbf{A} \bolds{\beta}(t)\bigr\} \,dt = O_p\bigl
\{d_n \bigl(\varphi_n \rho_n +
n^{-1/2} + d_nb_n^2\bigr)\bigr\}.
\]
Since $nb_n^{1/2} \varphi_n^2 \rho_n = o(1)$, (\ref{eqnPower})
follows from Lemma \ref{lemCLTTndiamond}.\vadjust{\goodbreak}
\end{pf*}

Recall that $D \subseteq D^*$ is a subset with complement $\bar
D$, and $D_0$ is the true set of relevant predictors. Let $\tilde
e_{D,i} = y_i - \mathbf{x}_{D,i}^\top\tilde{\bolds
\beta}(i/n)$, $1\le i \le n$. Then $\operatorname{RSS}(D) = \sum_{i=1}^n
\tilde e_{D,i}^2$. Lemma \ref{lemRSS} provides bounds for
$\operatorname{RSS}(D) - \sum_{i=1}^n e_i^2$ for both cases $D_0
\subseteq D$ and $D_0 \not\subseteq D$.

%
\begin{lemma}\label{lemRSS}
Assume \textup{(A1)--(A4)}, $\Theta_{0,4}(\mathbf{L}) < \infty$,
$\Theta_{n,\iota}(\mathbf{L}) = \mathcal{O}(n^{-\nu})$ for some
$\nu> 1/2 - 1/\iota$, $b_n \to0$ and $nb_n \to\infty$. Then \textup{(i)}
if $D_0 \subseteq D$, then
\[
\operatorname{RSS}(D) = \sum_{i=1}^n
e_i^2 + O_p\bigl\{n\varphi_n (
\varphi_n + \rho_n)\bigr\};
\]
and \textup{(ii)} if $D_0 \not\subseteq D$, then
\begin{eqnarray*}
\operatorname{RSS}(D) & = & \sum_{i=1}^n
e_i^2 + \sum_{i=1}^n
\bolds \beta_{\bar D}(i/n)^\top E\bigl\{\mathbf{x}_{\bar D,i}
\mathbf{x}_{\bar D,i}^\top\bigr\} \bolds{\beta}_{\bar D}(i/n)
\\
& &{} + O_p\bigl\{n^{1/2} + n\varphi_n (
\varphi_n + \rho_n)\bigr\}.
\end{eqnarray*}
\end{lemma}

\begin{pf}
For (i), since $D_0 \subseteq D$, we have $\tilde e_{D,i} = e_i +
\mathbf{x}_{D,i}^\top\{\bolds{\beta}_D(i/n) -\break
\tilde{\bolds{\beta}}_D(i/n)\}$ and
\[
\operatorname{RSS}(D) = \sum_{i=1}^n \tilde
e_{D,i}^2 = \sum_{i=1}^n
e_i^2 + I_n - 2\II_n,
\]
where, by the proof of Theorem \ref{thmCLTahat}, $I_n =
\sum_{i=1}^n [\mathbf{x}_{D,i}^\top\{\tilde{\bolds
\beta}_D(i/n) - \bolds{\beta}_D(i/n)\}]^2 = O_p(n\varphi_n^2)$
and, by (\ref{eqnBahadurrep}) and Lemmas \ref{lembndlin} and
\ref{lemsupbnd},
\begin{eqnarray*}
\II_n & = & \sum_{i=1}^n (
\mathbf{x}_i e_i)^\top\mathbf{A}_D^\top
\mathbf{A}_D \bigl\{\tilde{\bolds{\beta}}(i/n) - \bolds{\beta}(i/n)
\bigr\}
\\
& = & \sum_{i=1}^n (\mathbf{x}_i
e_i)^\top\mathbf{A}_D^\top
\mathbf{A}_D \mathbf{M}(\mathbf{G},i/n)^{-1} \mathbf
R_{n,0}(i/n)
\\
& &{} + O_p\bigl(nb_n \varphi_n + n
\varphi_n \rho_n + n^{1/2} b_n^2
\bigr).
\end{eqnarray*}
Since, by Lemma \ref{lembndquad},
\begin{eqnarray*}
& & \frac{1 }{ nb_n} \sum_{i=1}^n \sum
_{j=1}^n (\mathbf{x}_i
e_i)^\top\mathbf{A}_D^\top
\mathbf{A}_D \mathbf{M}(\mathbf{G},i/n)^{-1} (
\mathbf{x}_j e_j) K \biggl(\frac{j/n-i/n
}{ b_n} \biggr)
\\
&&\qquad =  O_p\bigl(b_n^{-1/2} + b_n^{-1}
\bigr),
\end{eqnarray*}
we have $\II_n = O_p(b_n^{-1} + n\varphi_n\rho_n + n^{1/2}b_n^2)$.
Since $b_n^{-1} + n^{1/2} b_n^2 = o\{n\varphi_n(\varphi_n +
\rho_n)\}$, (i) follows. For (ii), since $\tilde e_{D,i} = e_i +
\mathbf{x}_{D,i}^\top\{\bolds{\beta}_D(i/n) -
\tilde{\bolds{\beta}}_D(i/n)\} +\break \mathbf{x}_{\bar D,i}^\top
\bolds{\beta}_{\bar D}(i/n)$, we have
\[
\operatorname{RSS}(D) = \sum_{i=1}^n
e_i^2 + I_n^\circ+ 2
\II_n^\circ+ \III_n^\circ,
\]
where, by (i),
\[
I_n^\circ= \sum_{i=1}^n
\bigl[e_i + \mathbf{x}_{D,i}^\top \bigl\{\bolds{
\beta}_D(i/n) - \tilde{\bolds{\beta}}_D(i/n)\bigr\}
\bigr]^2 - \sum_{i=1}^n
e_i^2 = O_p\bigl\{n\varphi_n (
\varphi_n + \rho_n)\bigr\}
\]
and, by Lemma \ref{lembndlin} and the argument on the quantity
$\II_n$ in (i),
\begin{eqnarray*}
\II_n^\circ&=& \sum_{i=1}^n
\bigl[e_i + \mathbf{x}_{D,i}^\top \bigl\{\bolds{
\beta}_D(i/n) - \tilde{\bolds{\beta}}_D(i/n)\bigr\}
\bigr] \mathbf{x}_{\bar D,i}^\top\bolds{\beta}_{\bar D}(i/n)
\cr
&=& O_p\bigl(n^{1/2} + b_n^{-1} +
n\varphi_n\rho_n + n^{1/2}b_n^2
\bigr).
\end{eqnarray*}
In addition, by Lemma \ref{lembndlin},
\[
\III_n^\circ= \sum_{i=1}^n
\bigl\{\mathbf{x}_{\bar D,i}^\top \bolds{\beta}_{\bar D}(i/n)
\bigr\}^2 = \sum_{i=1}^n \bolds
\beta_{\bar D}(i/n)^\top E\bigl\{\mathbf{x}_{\bar D,i}
\mathbf{x}_{\bar D,i}^\top\bigr\} \bolds{\beta}_{\bar D}(i/n)
+ O_p\bigl(n^{1/2}\bigr),
\]
Lemma \ref{lemRSS} follows.
\end{pf}

\begin{pf*}{Proof of Theorem \ref{thmmodelselect}} By Lemma
\ref{lembndlin}, $\sum_{i=1}^n (e_i^2 - E e_i^2) = O_p(n^{1/2})$.
Lemma \ref{lemRSS} implies
\[
\log\bigl\{\operatorname{RSS}(D)\bigr\} = \log \Biggl(\sum
_{i=1}^n e_i^2 \Biggr) +
O_p\bigl\{\varphi_n(\varphi_n+
\rho_n)\bigr\}
\]
for $D_0 \subseteq D$, and
\[
\log\bigl\{\operatorname{RSS}(D)\bigr\} = \log \Biggl[\sum
_{i=1}^n e_i^2 + \sum
_{i=1}^n \bolds{\beta}_{\bar D}(i/n)^\top
E\bigl\{\mathbf{x}_{\bar D,i} \mathbf{x}_{\bar D,i}^\top\bigr
\} \bolds \beta_{\bar D}(i/n) \Biggr] + o_p(1)
\]
for $D_0 \not\subseteq D$. Since $\chi_n = o(1)$ and
$\varphi_n(\varphi_n + \rho_n) = o(\chi_n)$, Theorem
\ref{thmmodelselect} follows.
\end{pf*}

\begin{pf*}{Proof of Proposition \ref{thmcovest}}
By Lemma \ref{lembndlin},
\[
\sup_{t \in[0,1]} \bigl\|\hat{\mathbf{M}}(\mathbf{G}, t) - E\bigl\{\hat{
\mathbf{M}}(\mathbf{G}, t)\bigr\}\bigr\| = O\bigl\{(n\varpi_n)^{-1/2}
\bigr\},
\]
and, by Lemma \ref{lembndquad},
\[
\sup_{t \in[0,1]} \bigl\|\hat{\bolds{\Lambda}}(\mathbf{L}, t) - E\bigl\{\hat{
\bolds{\Lambda}}(\mathbf{L}, t)\bigr\}\bigr\| = O\bigl(\varrho_n^{1/2}
\bigr).
\]
By (\ref{eqnbndt1t2}), we have
\begin{eqnarray*}
& & \max_{1 \leq i \leq n} \bigl|E\bigl\{\bolds{\lambda}_i(\mathbf{L},
\tau_n \varrho_n)\bigr\} - \bolds{\Lambda}(\mathbf{L},
t_{i,n})\bigr|
\\
&&\qquad \leq C \sum_{k=0}^\infty\bigl\{
\gamma_k(\mathbf{L}) \wedge(\tau_n \varrho_n)
\bigr\} + \Theta_{n \tau_n \varrho_n,2}(\mathbf{L})
\\
&&\qquad =  O\bigl\{(\tau_n \varrho_n)^{\nu/(1+\nu)} + (n
\tau_n\varrho_n)^{-\nu}\bigr\}.
\end{eqnarray*}
Proposition \ref{thmcovest} follows by properties of local linear estimates.
\end{pf*}

\begin{pf*}{Proof of Proposition \ref{propTVAR}} Consider the process
$\{\mathbf z_{t,i}\}_{i \in\mathbb Z}$ that satisfies the
recursion
\[
\mathbf{z}_{t,i} = \mathbf{A}(t) \mathbf{z}_{t,i-1} +
\mathbf{H}^\diamond(t;\bolds{\mathcal F}_i),\qquad i \in \mathbb{Z}.
\]
Then, for each $t \in[0,1]$, the process $\{\mathbf
z_{t,i}\}_{i \in\mathbb{Z}}$ is stationary, and there exists a
measurable function $\mathbf{G}$ such that $\mathbf z_{t,i}
= \mathbf{G}(t; \bolds{\mathcal{F}}_i)$, $i \in
\mathbb{Z}$. By condition (T3), $\rho_{\mathbf{A}} = \sup_{t
\in[0,1]} \overline\rho\{\mathbf{A}(t)\} < 1$. Hence, by
condition (T2) and induction, we have
%
\begin{eqnarray}
\label{eqnTVARmaxdiffk} \max_{1 \leq i \leq n} \|\mathbf{x}_i -
\mathbf{z}_{i/n,i}\| &\leq&\rho_{\mathbf{A}}^k
\max_{1 \leq
i \leq n} \|\mathbf{x}_{i-k} - \mathbf{z}_{i/n,i-k}\|
\nonumber
\\[-8pt]
\\[-8pt]
\nonumber
&&{} +
C \sum_{j=1}^{k-1} \frac{j \rho_{\mathbf{A}}^j }{ n}, \qquad k
\geq 2.
\end{eqnarray}
Since $\rho_{\mathbf{A}} < 1$ and $\sum_{j=1}^\infty j
\rho_{\mathbf{A}}^j < \infty$, (\ref{eqnTVARmaxdiff}) follows
by letting $k \to\infty$. It suffices to show that $\mathbf{G}
\in \operatorname{Lip}$. For this, by a similar argument of
(\ref{eqnTVARmaxdiffk}), we have for any $k \geq2$,
\[
\sup_{t_1,t_2 \in[0,1]} \|\mathbf{z}_{t_1,i} - \mathbf{z}_{t_2,i}\|
\leq\rho_{\mathbf{A}}^k \sup_{t_1,t_2
\in[0,1]} \|
\mathbf{z}_{t_1,i-k} - \mathbf{z}_{t_2,i-k}\| + C |t_1-t_2|
\sum_{j=0}^{k-1} \rho_{\mathbf{A}}^j.
\]
Since $\rho_{\mathbf{A}} < 1$ and $\sum_{j=1}^\infty
\rho_{\mathbf{A}}^j < \infty$, Proposition \ref{propTVAR}
follows by letting $k \to\infty$.
\end{pf*}
\end{appendix}

\section*{Acknowledgments}
We are grateful to the Editor, an Associate Editor, and two anonymous
referees for their helpful comments and suggestions.

%


\printaddresses


\begin{thebibliography}{57}

\bibitem{AbramovichSpencerTurley2007}
%
\begin{barticle}[mr]
\bauthor{\bsnm{Abramovich},~\bfnm{Yuri~I.}\binits{Y.~I.}},
\bauthor{\bsnm{Spencer},~\bfnm{Nicholas~K.}\binits{N.~K.}} \AND
\bauthor{\bsnm{Turley},~\bfnm{Michael D.~E.}\binits{M.~D.~E.}}
(\byear{2007}).
\btitle{Order estimation and discrimination between stationary and time-varying
({TVAR}) autoregressive models}.
\bjournal{IEEE Trans. Signal Process.}
\bvolume{55}
\bpages{2861--2876}.
\bid{doi={10.1109/TSP.2007.893966}, issn={1053-587X}, mr={2473559}}
\bptok{imsref}%
\end{barticle}
%
\endbibitem

\bibitem{Andrews1993}
%
\begin{barticle}[mr]
\bauthor{\bsnm{Andrews},~\bfnm{Donald W.~K.}\binits{D.~W.~K.}}
(\byear{1993}).
\btitle{Tests for parameter instability and structural change with unknown
change point}.
\bjournal{Econometrica}
\bvolume{61}
\bpages{821--856}.
\bid{doi={10.2307/2951764}, issn={0012-9682}, mr={1231678}}
\bptok{imsref}%
\end{barticle}
%
\endbibitem

\bibitem{AndrewsMonahan1992}
%
\begin{barticle}[mr]
\bauthor{\bsnm{Andrews},~\bfnm{Donald W.~K.}\binits{D.~W.~K.}} \AND
\bauthor{\bsnm{Monahan},~\bfnm{J.~Christopher}\binits{J.~C.}}
(\byear{1992}).
\btitle{An improved heteroskedasticity and autocorrelation consistent
covariance matrix estimator}.
\bjournal{Econometrica}
\bvolume{60}
\bpages{953--966}.
\bid{doi={10.2307/2951574}, issn={0012-9682}, mr={1168742}}
\bptok{imsref}%
\end{barticle}
%
\endbibitem

\bibitem{BickelLevina2008}
%
\begin{barticle}[mr]
\bauthor{\bsnm{Bickel},~\bfnm{Peter~J.}\binits{P.~J.}} \AND
\bauthor{\bsnm{Levina},~\bfnm{Elizaveta}\binits{E.}}
(\byear{2008}).
\btitle{Regularized estimation of large covariance matrices}.
\bjournal{Ann. Statist.}
\bvolume{36}
\bpages{199--227}.
\bid{doi={10.1214/009053607000000758}, issn={0090-5364}, mr={2387969}}
\bptok{imsref}%
\end{barticle}
%
\endbibitem

\bibitem{BrownDurbinEvans1975}
%
\begin{barticle}[mr]
\bauthor{\bsnm{Brown},~\bfnm{R.~L.}\binits{R.~L.}},
\bauthor{\bsnm{Durbin},~\bfnm{James}\binits{J.}} \AND
\bauthor{\bsnm{Evans},~\bfnm{J.~M.}\binits{J.~M.}}
(\byear{1975}).
\btitle{Techniques for testing the constancy of regression
relationships over
time}.
\bjournal{J. Roy. Statist. Soc. Ser. B}
\bvolume{37}
\bpages{149--192}.
\bid{issn={0035-9246}, mr={0378310}}
\bptnote{check related}%
\bptok{imsref}%
\end{barticle}
%
\endbibitem

\bibitem{Cai2007}
%
\begin{barticle}[mr]
\bauthor{\bsnm{Cai},~\bfnm{Zongwu}\binits{Z.}}
(\byear{2007}).
\btitle{Trending time-varying coefficient time series models with serially
correlated errors}.
\bjournal{J. Econometrics}
\bvolume{136}
\bpages{163--188}.
\bid{doi={10.1016/j.jeconom.2005.08.004}, issn={0304-4076}, mr={2328589}}
\bptok{imsref}%
\end{barticle}
%
\endbibitem

\bibitem{CaiFanLi2000}
%
\begin{barticle}[mr]
\bauthor{\bsnm{Cai},~\bfnm{Zongwu}\binits{Z.}},
\bauthor{\bsnm{Fan},~\bfnm{Jianqing}\binits{J.}} \AND
\bauthor{\bsnm{Li},~\bfnm{Runze}\binits{R.}}
(\byear{2000}).
\btitle{Efficient estimation and inferences for varying-coefficient models}.
\bjournal{J. Amer. Statist. Assoc.}
\bvolume{95}
\bpages{888--902}.
\bid{doi={10.2307/2669472}, issn={0162-1459}, mr={1804446}}
\bptok{imsref}%
\end{barticle}
%
\endbibitem

\bibitem{ChenHong2012}
%
\begin{barticle}[auto:STB|2012/06/08|12:49:54]
\bauthor{\bsnm{Chen},~\bfnm{B.}\binits{B.}} \AND
\bauthor{\bsnm{Hong},~\bfnm{Y.}\binits{Y.}}
(\byear{2012}).
\btitle{Testing for smooth structural changes in time series
models via nonparametric regression}.
\bjournal{Econometrica}
\bvolume{80}
\bpages{1157--1183}.
\bptok{imsref}%
\end{barticle}
%
\endbibitem

\bibitem{Chow1960}
%
\begin{barticle}[mr]
\bauthor{\bsnm{Chow},~\bfnm{Gregory~C.}\binits{G.~C.}}
(\byear{1960}).
\btitle{Tests of equality between sets of coefficients in two linear
regressions}.
\bjournal{Econometrica}
\bvolume{28}
\bpages{591--605}.
\bid{issn={0012-9682}, mr={0141193}}
\bptok{imsref}%
\end{barticle}
%
\endbibitem

\bibitem{CravenWahba1979}
%
\begin{barticle}[mr]
\bauthor{\bsnm{Craven},~\bfnm{Peter}\binits{P.}} \AND
\bauthor{\bsnm{Wahba},~\bfnm{Grace}\binits{G.}}
(\byear{1979}).
\btitle{Smoothing noisy data with spline functions. {E}stimating the correct
degree of smoothing by the method of generalized cross-validation}.
\bjournal{Numer. Math.}
\bvolume{31}
\bpages{377--403}.
\bid{doi={10.1007/BF01404567}, issn={0029-599X}, mr={0516581}}
\bptnote{check year}%
\bptok{imsref}%
\end{barticle}
%
\endbibitem

\bibitem{Dahlhaus1996}
%
\begin{barticle}[mr]
\bauthor{\bsnm{Dahlhaus},~\bfnm{R.}\binits{R.}}
(\byear{1996}).
\btitle{On the {K}ullback--{L}eibler information divergence of locally
stationary processes}.
\bjournal{Stochastic Process. Appl.}
\bvolume{62}
\bpages{139--168}.
\bid{doi={10.1016/0304-4149(95)00090-9}, issn={0304-4149}, mr={1388767}}
\bptok{imsref}%
\end{barticle}
%
\endbibitem

\bibitem{Dahlhaus1997}
%
\begin{barticle}[mr]
\bauthor{\bsnm{Dahlhaus},~\bfnm{R.}\binits{R.}}
(\byear{1997}).
\btitle{Fitting time series models to nonstationary processes}.
\bjournal{Ann. Statist.}
\bvolume{25}
\bpages{1--37}.
\bid{doi={10.1214/aos/1034276620}, issn={0090-5364}, mr={1429916}}
\bptok{imsref}%
\end{barticle}
%
\endbibitem

\bibitem{DahlhausNeumannSachs1999}
%
\begin{barticle}[mr]
\bauthor{\bsnm{Dahlhaus},~\bfnm{Rainer}\binits{R.}},
\bauthor{\bsnm{Neumann},~\bfnm{Michael~H.}\binits{M.~H.}} \AND
\bauthor{\bparticle{von} \bsnm{Sachs},~\bfnm{Rainer}\binits{R.}}
(\byear{1999}).
\btitle{Nonlinear wavelet estimation of time-varying autoregressive processes}.
\bjournal{Bernoulli}
\bvolume{5}
\bpages{873--906}.
\bid{doi={10.2307/3318448}, issn={1350-7265}, mr={1715443}}
\bptok{imsref}%
\end{barticle}
%
\endbibitem

\bibitem{DavisHuangYao1995}
%
\begin{barticle}[mr]
\bauthor{\bsnm{Davis},~\bfnm{Richard~A.}\binits{R.~A.}},
\bauthor{\bsnm{Huang},~\bfnm{Da~Wei}\binits{D.~W.}} \AND
\bauthor{\bsnm{Yao},~\bfnm{Yi-Ching}\binits{Y.-C.}}
(\byear{1995}).
\btitle{Testing for a change in the parameter values and order of an
autoregressive model}.
\bjournal{Ann. Statist.}
\bvolume{23}
\bpages{282--304}.
\bid{doi={10.1214/aos/1176324468}, issn={0090-5364}, mr={1331669}}
\bptok{imsref}%
\end{barticle}
%
\endbibitem

\bibitem{DetteSpreckelsen2004}
%
\begin{barticle}[mr]
\bauthor{\bsnm{Dette},~\bfnm{Holger}\binits{H.}} \AND
\bauthor{\bsnm{Spreckelsen},~\bfnm{Ingrid}\binits{I.}}
(\byear{2004}).
\btitle{Some comments on specification tests in nonparametric absolutely
regular processes}.
\bjournal{J. Time Series Anal.}
\bvolume{25}
\bpages{159--172}.
\bid{doi={10.1111/j.1467-9892.2004.00343.x}, issn={0143-9782}, mr={2045571}}
\bptok{imsref}%
\end{barticle}
%
\endbibitem

\bibitem{FanHuang2005}
%
\begin{barticle}[mr]
\bauthor{\bsnm{Fan},~\bfnm{Jianqing}\binits{J.}} \AND
\bauthor{\bsnm{Huang},~\bfnm{Tao}\binits{T.}}
(\byear{2005}).
\btitle{Profile likelihood inferences on semiparametric varying-coefficient
partially linear models}.
\bjournal{Bernoulli}
\bvolume{11}
\bpages{1031--1057}.
\bid{doi={10.3150/bj/1137421639}, issn={1350-7265}, mr={2189080}}
\bptok{imsref}%
\end{barticle}
%
\endbibitem

\bibitem{FanYaoCai2003}
%
\begin{barticle}[mr]
\bauthor{\bsnm{Fan},~\bfnm{Jianqing}\binits{J.}},
\bauthor{\bsnm{Yao},~\bfnm{Qiwei}\binits{Q.}} \AND
\bauthor{\bsnm{Cai},~\bfnm{Zongwu}\binits{Z.}}
(\byear{2003}).
\btitle{Adaptive varying-coefficient linear models}.
\bjournal{J.~R.~Stat. Soc. Ser. B Stat. Methodol.}
\bvolume{65}
\bpages{57--80}.
\bid{doi={10.1111/1467-9868.00372}, issn={1369-7412}, mr={1959093}}
\bptok{imsref}%
\end{barticle}
%
\endbibitem

\bibitem{FanZhangZhang2001}
%
\begin{barticle}[mr]
\bauthor{\bsnm{Fan},~\bfnm{Jianqing}\binits{J.}},
\bauthor{\bsnm{Zhang},~\bfnm{Chunming}\binits{C.}} \AND
\bauthor{\bsnm{Zhang},~\bfnm{Jian}\binits{J.}}
(\byear{2001}).
\btitle{Generalized likelihood ratio statistics and {W}ilks phenomenon}.
\bjournal{Ann. Statist.}
\bvolume{29}
\bpages{153--193}.
\bid{doi={10.1214/aos/996986505}, issn={0090-5364}, mr={1833962}}
\bptok{imsref}%
\end{barticle}
%
\endbibitem

\bibitem{FanZhang1999}
%
\begin{barticle}[mr]
\bauthor{\bsnm{Fan},~\bfnm{Jianqing}\binits{J.}} \AND
\bauthor{\bsnm{Zhang},~\bfnm{Wenyang}\binits{W.}}
(\byear{1999}).
\btitle{Statistical estimation in varying coefficient models}.
\bjournal{Ann. Statist.}
\bvolume{27}
\bpages{1491--1518}.
\bid{doi={10.1214/aos/1017939139}, issn={0090-5364}, mr={1742497}}
\bptok{imsref}%
\end{barticle}
%
\endbibitem

\bibitem{FanZhang2000}
%
\begin{barticle}[mr]
\bauthor{\bsnm{Fan},~\bfnm{Jianqing}\binits{J.}} \AND
\bauthor{\bsnm{Zhang},~\bfnm{Wenyang}\binits{W.}}
(\byear{2000}).
\btitle{Simultaneous confidence bands and hypothesis testing in
varying-coefficient models}.
\bjournal{Scand. J. Stat.}
\bvolume{27}
\bpages{715--731}.
\bid{doi={10.1111/1467-9469.00218}, issn={0303-6898}, mr={1804172}}
\bptok{imsref}%
\end{barticle}
%
\endbibitem

\bibitem{FanLinton2003}
%
\begin{barticle}[mr]
\bauthor{\bsnm{Fan},~\bfnm{Yanqin}\binits{Y.}} \AND
\bauthor{\bsnm{Linton},~\bfnm{Oliver}\binits{O.}}
(\byear{2003}).
\btitle{Some higher-order theory for a consistent non-parametric model
specification test}.
\bjournal{J. Statist. Plann. Inference}
\bvolume{109}
\bpages{125--154}.
\bid{doi={10.1016/S0378-3758(02)00307-5}, issn={0378-3758}, mr={1946644}}
\bptok{imsref}%
\end{barticle}
%
\endbibitem

\bibitem{GaoGijbels2008}
%
\begin{barticle}[mr]
\bauthor{\bsnm{Gao},~\bfnm{Jiti}\binits{J.}} \AND
\bauthor{\bsnm{Gijbels},~\bfnm{Ir{\`e}ne}\binits{I.}}
(\byear{2008}).
\btitle{Bandwidth selection in nonparametric kernel testing}.
\bjournal{J. Amer. Statist. Assoc.}
\bvolume{103}
\bpages{1584--1594}.
\bid{doi={10.1198/016214508000000968}, issn={0162-1459}, mr={2504206}}
\bptok{imsref}%
\end{barticle}
%
\endbibitem

\bibitem{GaoHawthorne2006}
%
\begin{barticle}[mr]
\bauthor{\bsnm{Gao},~\bfnm{Jiti}\binits{J.}} \AND
\bauthor{\bsnm{Hawthorne},~\bfnm{Kim}\binits{K.}}
(\byear{2006}).
\btitle{Semiparametric estimation and testing of the trend of temperature
series}.
\bjournal{Econom. J.}
\bvolume{9}
\bpages{332--355}.
\bid{doi={10.1111/j.1368-423X.2006.00188.x}, issn={1368-4221}, mr={2324973}}
\bptok{imsref}%
\end{barticle}
%
\endbibitem

\bibitem{GencagaKuruogluErtuzunYildirim2008}
%
\begin{barticle}[auto:STB|2012/06/08|12:49:54]
\bauthor{\bsnm{Gen{\c{c}}a{\u{g}}a},~\bfnm{D.}\binits{D.}},
\bauthor{\bsnm{Kuruo{\u{g}}lu},~\bfnm{E.~E.}\binits{E.~E.}},
\bauthor{\bsnm{Ert{\"u}z{\"u}n},~\bfnm{A.}\binits{A.}} \AND
\bauthor{\bsnm{Yildirim},~\bfnm{S.}\binits{S.}}
(\byear{2008}).
\btitle{Estimation of time-varying $\operatorname{AR} S \alpha S$
processes using Gibbs
sampling}.
\bjournal{Signal Processing}
\bvolume{88}
\bpages{2564--2572}.
\bptok{imsref}%
\end{barticle}
%
\endbibitem

\bibitem{Grenier1983}
%
\begin{barticle}[auto:STB|2012/06/08|12:49:54]
\bauthor{\bsnm{Grenier},~\bfnm{Y.}\binits{Y.}}
(\byear{1983}).
\btitle{Time-dependent ARMA modelling of nonstationary signals}.
\bjournal{IEEE Trans. Acoust. Speech}
\bvolume{31}
\bpages{899--911}.
\bptok{imsref}%
\end{barticle}
%
\endbibitem

\bibitem{Hausman1978}
%
\begin{barticle}[mr]
\bauthor{\bsnm{Hausman},~\bfnm{J.~A.}\binits{J.~A.}}
(\byear{1978}).
\btitle{Specification tests in econometrics}.
\bjournal{Econometrica}
\bvolume{46}
\bpages{1251--1271}.
\bid{doi={10.2307/1913827}, issn={0012-9682}, mr={0513692}}
\bptok{imsref}%
\end{barticle}
%
\endbibitem

\bibitem{HeTerasvirtaGonzalez2009}
%
\begin{barticle}[mr]
\bauthor{\bsnm{He},~\bfnm{Changli}\binits{C.}},
\bauthor{\bsnm{Ter{\"a}svirta},~\bfnm{Timo}\binits{T.}} \AND
\bauthor{\bsnm{Gonz{\'a}lez},~\bfnm{Andr{\'e}s}\binits{A.}}
(\byear{2009}).
\btitle{Testing parameter constancy in stationary vector
autoregressive models
against continuous change}.
\bjournal{Econometric Rev.}
\bvolume{28}
\bpages{225--245}.
\bid{doi={10.1080/07474930802388041}, issn={0747-4938}, mr={2655626}}
\bptok{imsref}%
\end{barticle}
%
\endbibitem

\bibitem{HooverRiceWuYang1998}
%
\begin{barticle}[mr]
\bauthor{\bsnm{Hoover},~\bfnm{Donald~R.}\binits{D.~R.}},
\bauthor{\bsnm{Rice},~\bfnm{John~A.}\binits{J.~A.}},
\bauthor{\bsnm{Wu},~\bfnm{Colin~O.}\binits{C.~O.}} \AND
\bauthor{\bsnm{Yang},~\bfnm{Li-Ping}\binits{L.-P.}}
(\byear{1998}).
\btitle{Nonparametric smoothing estimates of time-varying coefficient models
with longitudinal data}.
\bjournal{Biometrika}
\bvolume{85}
\bpages{809--822}.
\bid{doi={10.1093/biomet/85.4.809}, issn={0006-3444}, mr={1666699}}
\bptok{imsref}%
\end{barticle}
%
\endbibitem

\bibitem{HuangWuZhou2004}
%
\begin{barticle}[mr]
\bauthor{\bsnm{Huang},~\bfnm{Jianhua~Z.}\binits{J.~Z.}},
\bauthor{\bsnm{Wu},~\bfnm{Colin~O.}\binits{C.~O.}} \AND
\bauthor{\bsnm{Zhou},~\bfnm{Lan}\binits{L.}}
(\byear{2004}).
\btitle{Polynomial spline estimation and inference for varying coefficient
models with longitudinal data}.
\bjournal{Statist. Sinica}
\bvolume{14}
\bpages{763--788}.
\bid{issn={1017-0405}, mr={2087972}}
\bptok{imsref}%
\end{barticle}
%
\endbibitem

\bibitem{KulasekeraWang1997}
%
\begin{barticle}[mr]
\bauthor{\bsnm{Kulasekera},~\bfnm{K.~B.}\binits{K.~B.}} \AND
\bauthor{\bsnm{Wang},~\bfnm{J.}\binits{J.}}
(\byear{1997}).
\btitle{Smoothing parameter selection for power optimality in testing of
regression curves}.
\bjournal{J. Amer. Statist. Assoc.}
\bvolume{92}
\bpages{500--511}.
\bid{doi={10.2307/2965699}, issn={0162-1459}, mr={1467844}}
\bptok{imsref}%
\end{barticle}
%
\endbibitem

\bibitem{LeybourneMcCabe1989}
%
\begin{barticle}[mr]
\bauthor{\bsnm{Leybourne},~\bfnm{S.~J.}\binits{S.~J.}} \AND
\bauthor{\bsnm{McCabe},~\bfnm{B.~P.~M.}\binits{B.~P.~M.}}
(\byear{1989}).
\btitle{On the distribution of some test statistics for coefficient constancy}.
\bjournal{Biometrika}
\bvolume{76}
\bpages{169--177}.
\bid{doi={10.1093/biomet/76.1.169}, issn={0006-3444}, mr={0991435}}
\bptok{imsref}%
\end{barticle}
%
\endbibitem

\bibitem{LinTerasvirta1994}
%
\begin{barticle}[auto:STB|2012/06/08|12:49:54]
\bauthor{\bsnm{Lin},~\bfnm{C.~F.~J.}\binits{C.~F.~J.}} \AND
\bauthor{\bsnm{Ter{\"a}svirta},~\bfnm{T.}\binits{T.}}
(\byear{1994}).
\btitle{Testing the constancy of regression parameters against continuous
structural change}.
\bjournal{J. Econometrics}
\bvolume{62}
\bpages{211--228}.
\bptok{imsref}%
\end{barticle}
%
\endbibitem

\bibitem{LinYing2001}
%
\begin{barticle}[mr]
\bauthor{\bsnm{Lin},~\bfnm{D.~Y.}\binits{D.~Y.}} \AND
\bauthor{\bsnm{Ying},~\bfnm{Z.}\binits{Z.}}
(\byear{2001}).
\btitle{Semiparametric and nonparametric regression analysis of longitudinal
data}.
\bjournal{J. Amer. Statist. Assoc.}
\bvolume{96}
\bpages{103--126}.
\bid{doi={10.1198/016214501750333018}, issn={0162-1459}, mr={1952726}}
\bptnote{check related}%
\bptok{imsref}%
\end{barticle}
%
\endbibitem

\bibitem{LiuWu2010}
%
\begin{barticle}[mr]
\bauthor{\bsnm{Liu},~\bfnm{Weidong}\binits{W.}} \AND
\bauthor{\bsnm{Wu},~\bfnm{Wei~Biao}\binits{W.~B.}}
(\byear{2010}).
\btitle{Asymptotics of spectral density estimates}.
\bjournal{Econometric Theory}
\bvolume{26}
\bpages{1218--1245}.
\bid{doi={10.1017/S026646660999051X}, issn={0266-4666}, mr={2660298}}
\bptok{imsref}%
\end{barticle}
%
\endbibitem

\bibitem{LiuXiaoWu2012}
%
\begin{bmisc}[auto:STB|2012/06/08|12:49:54]
\bauthor{\bsnm{Liu},~\bfnm{W.}\binits{W.}},
\bauthor{\bsnm{Xiao},~\bfnm{H.}\binits{H.}} \AND
\bauthor{\bsnm{Wu},~\bfnm{W.~B.}\binits{W.~B.}}
(\byear{2012}).
\bhowpublished{Probability and moment inequalities under dependence.
Working paper}.
\bptok{imsref}%
\end{bmisc}
%
\endbibitem

\bibitem{LumleyHeagerty1999}
%
\begin{barticle}[mr]
\bauthor{\bsnm{Lumley},~\bfnm{Thomas}\binits{T.}} \AND
\bauthor{\bsnm{Heagerty},~\bfnm{Patrick}\binits{P.}}
(\byear{1999}).
\btitle{Weighted empirical adaptive variance estimators for correlated data
regression}.
\bjournal{J. R. Stat. Soc. Ser. B Stat. Methodol.}
\bvolume{61}
\bpages{459--477}.
\bid{doi={10.1111/1467-9868.00187}, issn={1369-7412}, mr={1680302}}
\bptok{imsref}%
\end{barticle}
%
\endbibitem

\bibitem{MoulinesPriouretRoueff2005}
%
\begin{barticle}[mr]
\bauthor{\bsnm{Moulines},~\bfnm{Eric}\binits{E.}},
\bauthor{\bsnm{Priouret},~\bfnm{Pierre}\binits{P.}} \AND
\bauthor{\bsnm{Roueff},~\bfnm{Fran{\c{c}}ois}\binits{F.}}
(\byear{2005}).
\btitle{On recursive estimation for time varying autoregressive processes}.
\bjournal{Ann. Statist.}
\bvolume{33}
\bpages{2610--2654}.
\bid{doi={10.1214/009053605000000624}, issn={0090-5364}, mr={2253097}}
\bptok{imsref}%
\end{barticle}
%
\endbibitem

\bibitem{NabeyaTanaka1988}
%
\begin{barticle}[mr]
\bauthor{\bsnm{Nabeya},~\bfnm{Seiji}\binits{S.}} \AND
\bauthor{\bsnm{Tanaka},~\bfnm{Katsuto}\binits{K.}}
(\byear{1988}).
\btitle{Asymptotic theory of a test for the constancy of regression
coefficients against the random walk alternative}.
\bjournal{Ann. Statist.}
\bvolume{16}
\bpages{218--235}.
\bid{doi={10.1214/aos/1176350701}, issn={0090-5364}, mr={0924867}}
\bptok{imsref}%
\end{barticle}
%
\endbibitem

\bibitem{NeweyWest1987}
%
\begin{barticle}[mr]
\bauthor{\bsnm{Newey},~\bfnm{Whitney~K.}\binits{W.~K.}} \AND
\bauthor{\bsnm{West},~\bfnm{Kenneth~D.}\binits{K.~D.}}
(\byear{1987}).
\btitle{A simple, positive semidefinite, heteroskedasticity and autocorrelation
consistent covariance matrix}.
\bjournal{Econometrica}
\bvolume{55}
\bpages{703--708}.
\bid{doi={10.2307/1913610}, issn={0012-9682}, mr={0890864}}
\bptok{imsref}%
\end{barticle}
%
\endbibitem

\bibitem{Nyblom1989}
%
\begin{barticle}[mr]
\bauthor{\bsnm{Nyblom},~\bfnm{Jukka}\binits{J.}}
(\byear{1989}).
\btitle{Testing for the constancy of parameters over time}.
\bjournal{J. Amer. Statist. Assoc.}
\bvolume{84}
\bpages{223--230}.
\bid{issn={0162-1459}, mr={0999682}}
\bptok{imsref}%
\end{barticle}
%
\endbibitem

\bibitem{OrbeFerreiraRodriguez2005}
%
\begin{barticle}[mr]
\bauthor{\bsnm{Orbe},~\bfnm{Susan}\binits{S.}},
\bauthor{\bsnm{Ferreira},~\bfnm{Eva}\binits{E.}} \AND
\bauthor{\bsnm{Rodriguez-Poo},~\bfnm{Juan}\binits{J.}}
(\byear{2005}).
\btitle{Nonparametric estimation of time varying parameters under shape
restrictions}.
\bjournal{J. Econometrics}
\bvolume{126}
\bpages{53--77}.
\bid{doi={10.1016/j.jeconom.2004.02.006}, issn={0304-4076}, mr={2118278}}
\bptok{imsref}%
\end{barticle}
%
\endbibitem

\bibitem{PlobergerKramerKontrus1989}
%
\begin{barticle}[mr]
\bauthor{\bsnm{Ploberger},~\bfnm{Werner}\binits{W.}},
\bauthor{\bsnm{Kr{\"a}mer},~\bfnm{Walter}\binits{W.}} \AND
\bauthor{\bsnm{Kontrus},~\bfnm{Karl}\binits{K.}}
(\byear{1989}).
\btitle{A new test for structural stability in the linear regression model}.
\bjournal{J. Econometrics}
\bvolume{40}
\bpages{307--318}.
\bid{doi={10.1016/0304-4076(89)90087-0}, issn={0304-4076}, mr={0994952}}
\bptok{imsref}%
\end{barticle}
%
\endbibitem

\bibitem{RajanRaynerGodsill1997}
%
\begin{barticle}[auto:STB|2012/06/08|12:49:54]
\bauthor{\bsnm{Rajan},~\bfnm{J.~J.}\binits{J.~J.}},
\bauthor{\bsnm{Rayner},~\bfnm{P.~J.~W.}\binits{P.~J.~W.}} \AND
\bauthor{\bsnm{Godsill},~\bfnm{S.~J.}\binits{S.~J.}}
(\byear{1997}).
\btitle{Bayesian approach to parameter estimation and interpolation of
time-varying autoregressive processes using the Gibbs sampler}.
\bjournal{IEE P-Vis. Image Sign.}
\bvolume{144}
\bpages{249--256}.
\bptok{imsref}%
\end{barticle}
%
\endbibitem

\bibitem{RamsaySilverman2005}
%
\begin{bbook}[mr]
\bauthor{\bsnm{Ramsay},~\bfnm{J.~O.}\binits{J.~O.}} \AND
\bauthor{\bsnm{Silverman},~\bfnm{B.~W.}\binits{B.~W.}}
(\byear{2005}).
\btitle{Functional Data Analysis},
\bedition{2nd} ed.
\bpublisher{Springer}, \baddress{New York}.
\bid{mr={2168993}}
\bptok{imsref}%
\end{bbook}
%
\endbibitem

\bibitem{Robinson1989}
%
\begin{bincollection}[auto:STB|2012/06/08|12:49:54]
\bauthor{\bsnm{Robinson},~\bfnm{P.~M.}\binits{P.~M.}}
(\byear{1989}).
\btitle{Nonparametric estimation of time-varying parameters}.
In \bbooktitle{Statistical Analysis and Forecasting of Economic Structural
Change}
(\beditor{\bfnm{P.}\binits{P.}~\bsnm{Hackl}}, ed.)
\bpages{253--264}.
\bpublisher{Springer}, \baddress{Berlin}.
\bptok{imsref}%
\end{bincollection}
%
\endbibitem

\bibitem{Robinson1991}
%
\begin{bincollection}[auto:STB|2012/06/08|12:49:54]
\bauthor{\bsnm{Robinson},~\bfnm{P.~M.}\binits{P.~M.}}
(\byear{1991}).
\btitle{Time-varying nonlinear regression}.
In \bbooktitle{Economic Structure Change Analysis and Forecasting}
(\beditor{\bfnm{P.}\binits{P.}~\bsnm{Hackl}} \AND
\beditor{\bfnm{A.~H.}\binits{A.~H.}~\bsnm{Westland}}, eds.)
\bpages{179--190}.
\bpublisher{Springer}, \baddress{Berlin}.
\bptok{imsref}%
\end{bincollection}
%
\endbibitem

\bibitem{SubbaRao1970}
%
\begin{barticle}[mr]
\bauthor{\bsnm{Subba~Rao},~\bfnm{T.}\binits{T.}}
(\byear{1970}).
\btitle{The fitting of non-stationary time-series models with time-dependent
parameters}.
\bjournal{J. Roy. Statist. Soc. Ser. B}
\bvolume{32}
\bpages{312--322}.
\bid{issn={0035-9246}, mr={0269065}}
\bptok{imsref}%
\end{barticle}
%
\endbibitem

\bibitem{Wang2008}
%
\begin{barticle}[mr]
\bauthor{\bsnm{Wang},~\bfnm{Lan}\binits{L.}}
(\byear{2008}).
\btitle{Nonparametric test for checking lack of fit of the quantile regression
model under random censoring}.
\bjournal{Canad. J. Statist.}
\bvolume{36}
\bpages{321--336}.
\bid{doi={10.1002/cjs.5550360209}, issn={0319-5724}, mr={2431684}}
\bptok{imsref}%
\end{barticle}
%
\endbibitem

\bibitem{Wang1998}
%
\begin{barticle}[auto:STB|2012/06/08|12:49:54]
\bauthor{\bsnm{Wang},~\bfnm{Y.~D.}\binits{Y.~D.}}
(\byear{1998}).
\btitle{Smoothing spline models with correlated random errors}.
\bjournal{J.~Amer. Statist. Assoc.}
\bvolume{93}
\bpages{341--348}.
\bptok{imsref}%
\end{barticle}
%
\endbibitem

\bibitem{WuPourahmadi2009}
%
\begin{barticle}[mr]
\bauthor{\bsnm{Wu},~\bfnm{Wei~Biao}\binits{W.~B.}} \AND
\bauthor{\bsnm{Pourahmadi},~\bfnm{Mohsen}\binits{M.}}
(\byear{2009}).
\btitle{Banding sample autocovariance matrices of stationary processes}.
\bjournal{Statist. Sinica}
\bvolume{19}
\bpages{1755--1768}.
\bid{issn={1017-0405}, mr={2589209}}
\bptok{imsref}%
\end{barticle}
%
\endbibitem

\bibitem{XiaZhangTong2004}
%
\begin{barticle}[mr]
\bauthor{\bsnm{Xia},~\bfnm{Yingcun}\binits{Y.}},
\bauthor{\bsnm{Zhang},~\bfnm{Wenyang}\binits{W.}} \AND
\bauthor{\bsnm{Tong},~\bfnm{Howell}\binits{H.}}
(\byear{2004}).
\btitle{Efficient estimation for semivarying-coefficient models}.
\bjournal{Biometrika}
\bvolume{91}
\bpages{661--681}.
\bid{doi={10.1093/biomet/91.3.661}, issn={0006-3444}, mr={2090629}}
\bptok{imsref}%
\end{barticle}
%
\endbibitem

\bibitem{ZegerDiggle1994}
%
\begin{barticle}[pbm]
\bauthor{\bsnm{Zeger},~\bfnm{S.~L.}\binits{S.~L.}} \AND
\bauthor{\bsnm{Diggle},~\bfnm{P.~J.}\binits{P.~J.}}
(\byear{1994}).
\btitle{Semiparametric models for longitudinal data with application
to CD4
cell numbers in HIV seroconverters}.
\bjournal{Biometrics}
\bvolume{50}
\bpages{689--699}.
\bid{issn={0006-341X}, pmid={7981395}}
\bptok{imsref}%
\end{barticle}
%
\endbibitem

\bibitem{ZhangDette2004}
%
\begin{barticle}[mr]
\bauthor{\bsnm{Zhang},~\bfnm{Chunming}\binits{C.}} \AND
\bauthor{\bsnm{Dette},~\bfnm{Holger}\binits{H.}}
(\byear{2004}).
\btitle{A power comparison between nonparametric regression tests}.
\bjournal{Statist. Probab. Lett.}
\bvolume{66}
\bpages{289--301}.
\bid{doi={10.1016/j.spl.2003.11.005}, issn={0167-7152}, mr={2045474}}
\bptok{imsref}%
\end{barticle}
%
\endbibitem

\bibitem{ZhangWu2011}
%
\begin{barticle}[mr]
\bauthor{\bsnm{Zhang},~\bfnm{Ting}\binits{T.}} \AND
\bauthor{\bsnm{Wu},~\bfnm{Wei~Biao}\binits{W.~B.}}
(\byear{2011}).
\btitle{Testing parametric assumptions of trends of a~nonstationary time
series}.
\bjournal{Biometrika}
\bvolume{98}
\bpages{599--614}.
\bid{doi={10.1093/biomet/asr017}, issn={0006-3444}, mr={2836409}}
\bptok{imsref}%
\end{barticle}
%
\endbibitem

\bibitem{ZhangLeeSong2002}
%
\begin{barticle}[mr]
\bauthor{\bsnm{Zhang},~\bfnm{Wenyang}\binits{W.}},
\bauthor{\bsnm{Lee},~\bfnm{Sik-Yum}\binits{S.-Y.}} \AND
\bauthor{\bsnm{Song},~\bfnm{Xinyuan}\binits{X.}}
(\byear{2002}).
\btitle{Local polynomial fitting in semivarying coefficient model}.
\bjournal{J. Multivariate Anal.}
\bvolume{82}
\bpages{166--188}.
\bid{doi={10.1006/jmva.2001.2012}, issn={0047-259X}, mr={1918619}}
\bptok{imsref}%
\end{barticle}
%
\endbibitem

\bibitem{ZhouWu2010}
%
\begin{barticle}[mr]
\bauthor{\bsnm{Zhou},~\bfnm{Zhou}\binits{Z.}} \AND
\bauthor{\bsnm{Wu},~\bfnm{Wei~Biao}\binits{W.~B.}}
(\byear{2010}).
\btitle{Simultaneous inference of linear models with time varying
coefficients}.
\bjournal{J. R. Stat. Soc. Ser. B Stat. Methodol.}
\bvolume{72}
\bpages{513--531}.
\bid{doi={10.1111/j.1467-9868.2010.00743.x}, issn={1369-7412}, mr={2758526}}
\bptok{imsref}%
\end{barticle}
%
\endbibitem

\end{thebibliography}
\end{document}